\input amstex
\input xy
\xyoption{all}
\documentstyle{amsppt}
\document

\magnification 1000

\def\gen{\frak{g}}

\def\den{\frak{d}}
\def\ten{\frak{t}}

\def\een{\frak{e}}

\def\nen{\frak{n}}

\def\pen{\frak{p}}

\def\men{\frak{m}}

\def\Sen{\frak{S}}

\def\a{{\alpha}}
\def\g{{\gamma}}

\def\l{{\lambda}}
\def\b{{\beta}}

\def\ab{{\bold a}}
\def\bb{{\bold b}}

\def\fb{{\bold f}}

\def\ib{{\bold i}}
\def\ibt{{{\yb}}}

\def\kb{{\bold k}}

\def\mb{{\bold m}}

\def\yb{{\bold y}}

\def\Db{{\bold D}}

\def\Fb{{\bold F}}

\def\modb{{\pmb{\Mc od}}}
\def\shb{{\pmb{\Sc h}}}
\def\projb{{\pmb{\Pc roj}}}
\def\Dcb{{\pmb{\Dc}}}
\def\Pcb{{\pmb{\Pc}}}
\def\Qcb{{\pmb{\Qc}}}

\def\Kb{{\bold K}}

\def\Pb{{\bold P}}
\def\Qb{{\bold Q}}
\def\Rb{{\bold R}}

\def\Sb{{\bold S}}

\def\Vb{{\bold V}}

\def\Wb{{\bold W}}
\def\Yb{{\bold Y}}
\def\Zb{{\bold Z}}

\def\k{{\kb}}

\def\Ext{\roman{Ext}}

\def\Hom{{\roman{Hom}}}

\def\Eu{{\roman{eu}}}

\def\dim{{\roman{dim}}}

\def\Ind{{\roman{Ind}}}

\def\End{{\roman{End}}}

\def\eu{{\roman{eu}}}
\def\ch{{\roman{ch}}}

\def\Prod{\ts\prod}

\def\GL{{\roman{GL}}}

\def\CC{{\Bbb C}}

\def\NN{{\Bbb N}}

\def\PP{{\Bbb P}}
\def\QQ{{\Bbb Q}}

\def\ZZ{{\Bbb Z}}

\def\Ac{{\Cal A}}
\def\Bc{{\Cal B}}
\def\Cc{{\Cal C}}
\def\Dc{{\Cal D}}

\def\Hc{{\Cal H}}

\def\Kc{{\Cal K}}
\def\Lc{{\Cal L}}
\def\Mc{{\Cal M}}

\def\Oc{{\Cal O}}
\def\Pc{{\Cal P}}
\def\Qc{{\Cal Q}}

\def\Sc{{\Cal S}}

\def\Vc{{\Cal V}}
\def\Vcb{{\pmb{\Vc}}}

\def\Ccb{{\pmb{\Cc}}}

\def\and{{\text{and}}}

\def\ts{\textstyle}

\def\qed{\hfill $\sqcap \hskip-6.5pt \sqcup$}        
\overfullrule=0pt                                    

\def\7dag{{{\!\!\!\!\!\!\!\dag}}}
\def\6dag{{{\!\!\!\!\!\!\dag}}}
\def\5dag{{{\!\!\!\!\!\dag}}}
\def\4dag{{{\!\!\!\!\dag}}}
\def\3dag{{{\!\!\!\dag}}}
\def\2dag{{{\!\!\dag}}}
\def\1dag{{{\!\dag}}}

\def\ind{{\lim\limits_{\lra}}}

\def\lra{{{\longrightarrow}}}

\newdimen\Squaresize\Squaresize=14pt
\newdimen\Thickness\Thickness=0.5pt
\def\Square#1{\hbox{\vrule width\Thickness
          \alphaox to \Squaresize{\hrule height \Thickness\vss
          \hbox to \Squaresize{\hss#1\hss}
          \vss\hrule height\Thickness}
          \unskip\vrule width \Thickness}
          \kern-\Thickness}
\def\Vsquare#1{\alphaox{\Square{$#1$}}\kern-\Thickness}

\nologo

\topmatter
\title Canonical bases and KLR-algebras
\endtitle
\rightheadtext{} \abstract We prove a recent conjecture of
Khovanov-Lauda concerning the categorification of one-half of the
quantum group associated with a simply laced Cartan datum.
\endabstract
\author M. Varagnolo, E. Vasserot\endauthor
\address D\'epartement de Math\'ematiques,
Universit\'e de Cergy-Pontoise, 2 av. A. Chauvin, BP 222, 95302
Cergy-Pontoise Cedex, France, Fax : 01 34 25 66 45\endaddress \email
michela.varagnolo\@math.u-cergy.fr\endemail
\address D\'epartement de Math\'ematiques,
Universit\'e Paris 7, 175 rue du Chevaleret, 75013 Paris, France,
Fax : 01 44 27 78 18
\endaddress
\email vasserot\@math.jussieu.fr\endemail
\thanks
2000{\it Mathematics Subject Classification.} Primary ??; Secondary
??.
\endthanks
\endtopmatter
\document

\head \endhead

\subhead 0. Introduction and notation\endsubhead

In \cite{KL1}, \cite{KL2} Khovanov and Lauda have introduced a new
family of algebras and formulated some conjecture predicting a
connection between the representation theory of these algebras and
Lusztig's geometric construction of canonical bases. The goal of
this paper is to prove part of this conjecture.

The construction in loc.~cit.~is as follows. Let ${}_\Ac\fb$ be
Lusztig's integral form of the negative half of the quantum
universal enveloping algebra associated with a quiver $\Gamma$. Let
$I$ be the set of vertices of $\Gamma$. One construct a family of
graded rings $\Rb_\nu$ over $\nu\in\NN I$. These rings are defined
in a combinatorial way. One consider the  Grothendieck group
$K(\Rb_\nu)$ of the category of finitely generated graded projective
$\Rb_\nu$-modules. Set $\Ac=\ZZ[q,q^{-1}]$. The direct sum
$$K(\Rb)=\bigoplus_{\nu\in\NN I}K(\Rb_\nu)$$
has a natural structure of a free $\Ac$-module. It is a
$\Ac$-bialgebra under induction and restriction. In loc.~cit.~ an
$\Ac$-algebra isomorphism $\g_\Ac:{}_\Ac\fb\to K(\Rb)$ is given. It
is conjectured that $\g_\Ac^{-1}$ maps the classes of the indecomposable
projective modules to the canonical basis. We prove this here.

Before to go on let us make a few historical remarks. The
KLR-algebras were first introduced by Khovanov and Lauda in
\cite{KL1}, \cite{KL2} with some restrictions on the quiver. There
were independently discovered by Rouquier, in a more general
version. See Remark 3.3 below for details and \cite{R} for the
definition and the first properties. After our paper was written
Rouquier informed us he has also obtained the same result as ours,
independently.

Now we give some notation we'll use in this paper. By a canonical
isomorphism in a given category we'll mean an explicit isomorphism.
We'll identify two objects with a canonical isomorphism whenever
needed. Unless specified otherwise, by a graded object of an
additive category $\Ccb$ we'll always mean a $\ZZ$-graded object,
i.e., an object of the form $M=\bigoplus_{i\in\ZZ}M_i$ where each
$M_i$ is an object of $\Ccb$. Then, given an integer $r$, we'll
write $M[r]$ for the shift of the grading on $M$ up by $r$, i.e.,
the degree $i$ component of $M[r]$ is equal to $M_{i+r}$.

Given a positive integer $m$ and a tuple $\mb=(m_1,m_2,\dots m_r)$
of positive integers we write $\Sen_m$ for the symmetric group and
$\Sen_\mb$ for the group $\prod_{l=1}^r\Sen_{m_l}$. Set
$$|\mb|=\sum_{l=1}^rm_l,\quad
\ell_\mb=\sum_{l=1}^r\ell_{m_l},\quad\ell_m=m(m-1)/2.$$ We use the
following notation for $q$-numbers
$$[m]_q=\sum_{l=1}^mq^{m+1-2l},
\quad [m]_q!=\prod_{l=1}^m[l]_q, \quad
[\mb]_q!=\prod_{l=1}^r[m_l]_q!.$$ If no confusion is possible we'll
abbreviate
$$[m]=[m]_q, \quad [m]!= [m]_q!,
\quad [\mb]!=[\mb]_q!.$$ Given two tuples $\mb=(m_1,m_2,\dots m_r)$,
$\mb'=(m'_1,m'_2,\dots m'_{r'})$ we define the tuple
$$\mb\mb'=(m_1,m_2,\dots m_r,m'_1,m'_2,\dots m'_{r'}).\leqno(0.1)$$

Let $[\Ccb]$ denote the Grothendieck group of an exact category
$\Ccb$. Let $\kb$ be a field and $\Rb=\bigoplus_i\Rb_i$ be a
$\ZZ$-graded $\kb$-algebra. Let $\Rb$-$\modb$ be the
category of finitely generated graded left $\Rb$-modules,
$\Rb$-$\modb^f$ be the full subcategory of finite-dimensional modules and
$\Rb$-$\projb$ be the full subcategory of projective objects. We'll
abbreviate
$$K(\Rb)=[\Rb\text{-}\projb],\quad
G(\Rb)=[\Rb\text{-}\modb^f].$$ Assume that the $\kb$-vector spaces
$\Rb_i$ are finite dimensional for each $i$.
Then $K(\Rb)$ is a
free Abelian group with a basis formed by the isomorphism classes of
the indecomposable objects in $\Rb$-$\projb$. 
Further $G(\Rb)$ is a free Abelian group with a basis
formed by the isomorphism classes of the simple objects in
$\Rb$-$\modb^f$. Given an object $M$ of $\Rb$-$\projb$ or
$\Rb$-$\modb^f$ let $[M]$ denote its class in $K(\Rb)$, $G(\Rb)$
respectively.
Both $K(\Rb)$, $G(\Rb)$ are $\Ac$-modules, where $\Ac$ acts on
$G(\Rb)$, $K(\Rb)$ as follows
$$q[M]=[M[1]],\quad q^{-1}[M]=[M[-1]],\quad\forall M.\leqno(0.2)$$
If the grading of $\Rb$ is bounded below then the $\Ac$-modules
$K(\Rb)$, $G(\Rb)$ are  free.

From now on the symbol $\k$ will denote a field of
characteristic zero.

\subhead 1. Reminder on quivers, extensions and convolution algebras
\endsubhead

We give a few notation on quivers and equivariant homology. Given a
quiver $\Gamma$ we first recall the definition of the semisimple
complexes on the moduli stack of representations of $\Gamma$ introduced by
Lusztig in \cite{L1}. Next we define their Ext-algebras (with respect to the
Yoneda product). Finally we recall a lemma of Ginzburg which relates
Ext-algebras to convolution algebras in equivariant homology.

\subhead{1.1. Representations of quivers}\endsubhead 
Let $\Gamma$ be a finite nonempty quiver
such that no arrow may join a vertex to itself. Recall that
$\Gamma$ is a tuple $(I,H,h\mapsto h',h\mapsto h'')$ where $I$ is
the set of vertices, $H$ is the set of arrows and for each $h\in H$
the vertices $h',h''\in I$ are the origin and the goal of $h$
respectively. Although this hypothesis is not necessary, 
we'll always assume that the set $I$ is finite in the rest of the paper. 
For each $i,j\in I$ we write $$H_{i,j}=\{h\in
H;h'=i,h''=j\}.$$ We'll abbreviate $i\to j$ for
$H_{i,j}\neq\emptyset$, $i\not\to j$ for $H_{i,j}=\emptyset$, and
$h:i\to j$ for $h\in H_{i,j}$. Let $h_{i,j}$ be the number of
elements in $H_{i,j}$ and set
$$i\cdot j=-h_{i,j}-h_{j,i},\quad i\cdot i=2,\quad i\neq j.$$

Let $\Vcb$ be the category of finite-dimensional $I$-graded
$\CC$-vector spaces $\Vb=\bigoplus_{i\in I}\Vb_i$ with morphisms
being linear maps respecting the grading. For each $\nu=\sum_i\nu_i
i\in\NN I$ let $\Vcb_\nu$ be the full subcategory of $\Vcb$ whose
objects are those $\Vb$ such that $\dim(\Vb_i)=\nu_i$ for all $i$. We
call $\nu$ the dimension vector of $\Vb$.
Given $\Vb\in\Vcb$ let
$$E_\Vb=\bigoplus_{h\in H}\Hom(\Vb_{h'},\Vb_{h''}).$$
The algebraic group $G_\Vb=\Prod_i\GL(\Vb_i)$ acts on $E_\Vb$ by
$(g,x)\mapsto y$ where $y_h=g_{h''}x_hg_{h'}^{-1}$, $g=(g_i)$,
$x=(x_h)$, and $y=(y_h)$.

Fix a nonzero element $\nu$ of $\NN I$.
Let $Y_\nu$ be the set of all pairs $\ibt=(\ib,\ab)$ where
$\ib=(i_1,i_2,\dots i_k)$ is a sequence of elements of $I$ and
$\ab=(a_1,a_2,\dots a_k)$ is a sequence of positive integers such
that $\sum_la_l\,i_l=\nu$. The set $Y_\nu$ is finite. 
Let $I^\nu\subset Y_\nu$
be the set of all pairs $\ibt=(\ib,\ab)$ such that $a_l=1$ for all
$l$. We'll abbreviate $\ib$ for $\ibt=(\ib,\ab)$ if $\ibt\in I^\nu$.
Note that the assignment
$$\ibt\mapsto(a_1 i_1,a_2 i_2,\dots a_k i_k)$$
identifies $Y_\nu$ with a set of sequences
$$\nu^1,\nu^2,\dots,\nu^k\in\NN I\quad\roman{with}\ \nu=\sum_{l=1}^k\nu^l.
\leqno(1.1)$$
For each $\nu=\sum_{i\in I}\nu_i i\in\NN I$ we set
$|\nu|=\sum_i\nu_i$, an integer $\geqslant 0$. Given $m\in\NN$ we
have $\bigsqcup_{\nu}I^\nu=I^m,$ where $\nu$ runs over the subset
$\NN_m I\subset\NN I$ of elements with $|\nu|=m$.

From now on we assume that $\Vb\in\Vcb_\nu$. For any sequence 
$\yb=(\nu^1,\nu^2,\dots\nu^k)$ as in (1.1), a flag of type $\yb$ in
$\Vb$ is a sequence
$$\phi=(\Vb=\Vb^k\supset\Vb^{k-1}\supset\cdots\supset\Vb^0=0)$$ of
$I$-graded subspace of $\Vb$ such that for any $l$ the $I$-graded
subspace $\Vb^{l}/\Vb^{l-1}$ belongs to $\Vcb_{\nu^l}$. Let
$F_{\yb}$ be the variety of all flags of type $\yb$ in $\Vb$. The
group $G_\Vb$ acts transitively on $F_{\yb}$ in the obvious way,
yielding a smooth projective $G_\Vb$-variety structure on $F_\yb$.

We'll say that the sequence $\yb=(\nu^1,\nu^2,\dots\nu^k)$ as in (1.1)
is multiplicity free if we have $\nu^l=\sum_{i\in I}\nu^l_i i$ with
$\nu^l_i=0$ or $1$ for each $l$. Note that $I^\nu$ consists of
multiplicity free sequences. Given a finite dimensional
$\CC$-vector space $\Wb$ let $F(\Wb)$ be the variety of all complete
flags in $\Wb$. If $\yb$ is a multiplicity free sequence 
then we have
$$\dim(F_\yb)=\ell_\nu=\sum_{i\in I}\ell_{\nu_i},\quad
F_\yb\simeq\prod_{i\in I}F(\Vb_i),\quad
\phi\mapsto(\phi\cap\Vb_i).$$

If $x\in E_\Vb$ we say that the flag $\phi$ is $x$-stable if
$x_h(\Vb^l_{h'})\subset\Vb^l_{h''}$ for all $h$, $l$. Let
$\widetilde F_{\yb}$ be the variety of all pairs $(x,\phi)$ such
that $\phi$ is $x$-stable. Set
$$d_{\yb}=\dim(\widetilde F_{\yb}).\leqno(1.2)$$
The group $G_\Vb$ acts on $\widetilde F_{\yb}$ by
$g:(x,\phi)\mapsto(gx,g\phi)$. The first projection
$$\pi_{\yb}:\widetilde F_{\yb}\to E_\Vb$$
is a $G_\Vb$-equivariant proper morphism.

For each $l=1,2,\dots |\nu|$ we define $\Oc_{\widetilde F_\yb}(l)$
to be the $G_\Vb$-equivariant line bundle over $\widetilde F_\yb$
whose fiber at the flag $\phi$ is equal to $\Vb^{l}/\Vb^{l-1}$.

\subhead{1.2. Constructible sheaves}\endsubhead 
Given an action of a complex linear
algebraic group $G$ on a quasiprojective algebraic variety $X$ over
$\CC$ we write $\Dcb_G(X)$ for the bounded $G$-equivariant derived category of
complexes of sheaves of $\k$-vector spaces on $X$.
Objects of $\Dcb_G(X)$ are referred to as complexes. If $G=\{e\}$,
the trivial group, we abbreviate $\Dcb(X)=\Dcb_G(X)$. For each
complexes $\Lc$, $\Lc'$ we'll abbreviate
$$\Ext_G^*(\Lc,\Lc')=\Ext^*_{\Dcb_G(X)}(\Lc,\Lc'),\quad
\Ext^*(\Lc,\Lc')=\Ext^*_{\Dcb(X)}(\Lc,\Lc')$$ if no confusion is
possible.

There is a lot of literature on the category $\Dcb_G(X)$, see
\cite{BL}, \cite{J}, \cite{L3, sec.~1}, \cite{L4, sec.~1} for
instance. Although we'll only use standards properties of $\Dcb_G(X)$
let us recall a few definitions for the comfort of the reader. We'll
use the notation of \cite{BBD} for sheaves and Deligne's definition
of $\Dcb_G(X)$, see \cite{BL, sec.~I.2, app.~B}, \cite{J,
sec.~1.2.3}. More precisely, let $[G\setminus X]$ be the usual
simplicial topological set whose $n$-skeleton is $G^n\times X$ for
each $n$. We put the transcendental topology on $X$, $G$. Let
$\shb_G(X)$ be the full subcategory of the category of simplicial
sheaves on $[G\setminus X]$ for which all structure morphisms are
isomorphisms. It is an Abelian category which is equivalent to the
category of $G$-equivariant sheaves on $X$ by \cite{D, (6.1.2)}. We
define $\Dcb_G(X)$ as the full subcategory of the bounded derived
category of simplicial sheaves on $[X\setminus G]$ consisting of the
complexes whose cohomology sheaves belong to $\shb_G(X)$.

Let us recall a few simple facts. If $H\subset G$ is a closed
subgroup then the pull back by the obvious morphism of simplicial
topological spaces $[H\setminus X]\to[G\setminus X]$ is the
forgetful functor $\Dcb_G(X)\to \Dcb_{H}(X)$. An object of
$\Dcb_G(X)$ is perverse if the corresponding object in $\Dcb(X)$ is
perverse.


The constant sheaf $\k$ on $X$ will be denoted $\k$. For any object
$\Lc$ of $\Dcb_G(X)$ let $H_G^*(X,\Lc)$ be the space of
$G$-equivariant cohomology with coefficients in $\Lc$. Let $\Dc\in
\Dcb_G(X)$ be the $G$-equivariant dualizing complex, see \cite{BL,
def.~3.5.1}. For each $\Lc$ we'll abbreviate
$$\Lc^\vee=\Hc om(\Lc,\Dc)$$
(the internal Hom). Recall that
$$(\Lc^\vee)^\vee=\Lc,\quad
\Ext^*_G(\Lc,\Dc)=H^*_G(X,\Lc^\vee),\quad
\Ext^*_{G}(\k,\Lc)=H_G^*(X,\Lc).\leqno(1.3)$$ We define the space of
$G$-equivariant homology by
$$H^G_*(X,\k)=H^*_G(X,\Dc).\leqno(1.4)$$ Note that we have not chosen the
grading in the most usual way, compare \cite{L4, sec.~1.17} for
instance. Note also that $\Dc=\k[2d]$ if $X$ is a smooth $G$-variety
of pure dimension $d$. Consider the following graded $\kb$-algebra
$$\Sb_G=H^{*}_G(\bullet,\k).$$
The graded $\k$-vector space $H^G_*(X,\k)$ has a natural structure
of a graded $\Sb_G$-module. We have
$$H_*^G(\bullet,\k)=\Sb_G$$ as graded
$\Sb_G$-module. There is a canonical graded $\k$-algebra isomorphism
$$\Sb_G\simeq\k[\gen]^G.$$ Here the symbol $\gen$ denotes the Lie
algebra of $G$ and a $G$-invariant homogeneous polynomial over
$\gen$ of degree $d$ is given the degree $2d$ in $\Sb_G$.

Fix a morphism of quasi-projective algebraic $G$-varieties $f:X\to
Y$. If $f$ is a proper map there is a direct image homomorphism
$$f_*:H_*^G(X,\k)\to H_*^G(Y,\k).$$ If $f$ is a smooth map of relative
dimension $d$ there is an inverse image homomorphism
$$f^*:H_i^G(Y,\k)\to H_{i-2d}^G(X,\k),\quad\forall i.$$
If $X$ has pure dimension $d$ there is a natural homomorphism
$$H_G^i(X,\k)\to H^G_{i-2d}(X,\k).$$ It is invertible if $X$ is
smooth. The image of the unit is called the fundamental class of $X$
in $H_*^G(X,\k)$. We denote it by $[X]$. If $f:X\to Y$ is the
embedding of a $G$-stable closed subset and $X'\subset X$ is the
union of the irreducible components of maximal dimension then the
image of $[X']$ by the map $f_*$ is the fundamental class of $X$ in
$H_*^G(Y,\k)$. It is again denoted by $[X]$.

\subhead{1.3. Ext-algebras}\endsubhead 
In the rest of this section we assume that
$\nu\in\NN I$ and $\Vb\in\Vcb_\nu$. We abbreviate
$$\Sb_\Vb=\Sb_{G_\Vb}.$$ For each sequence $\yb$ as above we have
the following semisimple complexes in $\Dcb_{G_\Vb}(E_\Vb)$
$$\Lc_{\yb}=(\pi_{\yb})_!(\k),\quad\Lc_\yb^\vee=\Lc_\yb[2d_\yb],\quad
{}^\delta\!\Lc_\yb=\Lc_\yb[d_\yb].\leqno(1.5)$$ Here the integer
$d_\yb$ is as in (1.2). Given $\yb,\yb'\in Y_\nu$ we consider the
graded $\Sb_\Vb$-module
$$\Zb_{\yb,\yb'}=
\Ext^*_{G_\Vb}(\Lc_\yb^\vee,\Lc_{\yb'}^\vee).$$ Given $\yb''\in
Y_\nu$ the Yoneda composition is a homogeneous $\Sb_\Vb$-bilinear
map of degree zero
$$\star:\Zb_{\yb,\yb'}\times\Zb_{\yb',\yb''}\to\Zb_{\yb,\yb''}.$$
The map $\star$ equips the $\k$-vector spaces
$$\Zb'_\Vb=\bigoplus_{\ibt,\ibt'\in Y_\nu}
\Zb_{\ibt,\ibt'},\quad \Zb_\Vb=\bigoplus_{\ib,\ib'\in
I^\nu}\Zb_{\ib,\ib'}$$ with the structure of unital associative
graded $\Sb_\Vb$-algebras. Our first goal is to compute $\Zb'_\Vb$,
because it is the most relevant for Lusztig's construction of
canonical bases. Due to the following basic fact it is indeed enough
to consider only $\Zb_\Vb$, which is smaller.

\proclaim{1.4. Proposition} The categories $\Zb'_\Vb$-$\modb$ and
$\Zb_\Vb$-$\modb$ are equivalent.
\endproclaim

\noindent{\sl Proof :} Let $\ibt=(\ib,\ab)\in Y_\nu$. We have
$a_l>0$ for all $l$. If $a_l>1$ for some $l$, let $\ibt'\in Y_\nu$
be obtained from $\ibt$ by replacing the single entry $i_l$ by $a_l$
entries equal to $i_l$ and the single entry $a_l$ by $a_l$ entries
equal to 1. By \cite{L1, sec.~2.4} we have
$\Lc_{\ibt'}\simeq\sum_{l=1}^r\Lc_{\ibt}[d_l]$ for some integer
$d_1,d_2,\dots d_r$ with $r\geqslant 1$. Let us do this
simultaneously for all $l$ such that $a_l>1$. Let $\ib'\in I^\nu$ be
the sequence obtained by expanding the pair $\ibt$. So we have
$$\ib'=(i_1,\dots i_1,i_2,\dots i_2,\dots, i_k,\dots,i_k),\quad
\ib=(i_1,\dots,i_k),\quad \ab=(a_1,\dots, a_k),\leqno(1.6)$$ where
the entry $i_l$ have the multiplicity $a_l$ in $\ib'$. There are
integers $\ell_1,\ell_2,\dots,\ell_s$, with $s\geqslant 1$, such
that $\Lc_{\ib'}\simeq\sum_{l=1}^s\Lc_{\ibt}[\ell_l]$. Thus the
algebras $\Zb'_\Vb$, $\Zb_\Vb$ are Morita equivalent.

\qed

\vskip3mm

\subhead 1.5.~Remark\endsubhead The integers $\ell_1,\dots\ell_s$ in
the proof of Proposition 1.4 can be computed explicitly. Recall that
$\ibt=(\ib,\ab)\in  Y_\nu$ with $\ib=(i_1,\dots,i_k)$,
$\ab=(a_1,\dots, a_k)$ and that $\ib'\in I^\nu$ is the sequence
obtained by expanding the pair $\ibt$. In the category
$\Dcb_{G_\Vb}(E_\Vb)$ we have an isomorphism of complexes
$$\Lc_{\ib'}\simeq\bigoplus_{w\in\Sen_\ab}\Lc_{\ibt}[-2\ell(w)].$$
This isomorphism is not canonical. It depends of the choice of a
partition of the variety $\widetilde F_{\ib'}$ into cells as in
\cite{L2, sec.~8.1.6}.

\subhead 1.6.~The canonical module over $\Zb_\Vb$\endsubhead 
For each $\yb\in Y_\nu$ we have the graded
$\Sb_\Vb$-modules
$$\Fb_\yb=\Ext^*_{G_\Vb}(\Lc_\yb^\vee,\Dc), \quad
\Fb_\Vb=\bigoplus_{\ib\in I^\nu}\Fb_\ib.\leqno(1.7)$$  The Yoneda
product gives a graded $\Sb_\Vb$-bilinear map
$$\Zb_{\yb,\yb'}\times\Fb_{\yb'}\to\Fb_{\yb}.$$
This yields a left graded representation of $\Zb_\Vb$ on $\Fb_\Vb$.
For each $\ib\in I^\nu$ let $1_\ib\in\Zb_{\ib,\ib}$ denote the
identity of $\Lc_\ib$. The elements $1_\ib$ form a complete set of
orthogonal idempotents of $\Zb_\Vb$ such that
$$\Zb_{\ib,\ib'}=1_{\ib}\star\Zb_\Vb\star 1_{\ib'},\quad
\Fb_{\ib}=1_\ib\star\Fb_\Vb.$$

\subhead 1.7.~Convolution algebras\endsubhead 
For $\ib,\ib'\in I^\nu$ we set
$$Z_{\ib,\ib'}=\widetilde F_\ib\times_{E_\Vb}\widetilde F_{\ib'},$$
the reduced fiber product relative to the maps $\pi_\ib$,
$\pi_{\ib'}$. We set also $$Z_\Vb=\coprod_{\ib,\ib'\in
I^\nu}Z_{\ib,\ib'},\quad\widetilde F_\Vb=\coprod_{\ib\in
I^\nu}\widetilde F_\ib,\quad F_\Vb=\coprod_{\ib\in I^\nu} F_\ib.$$
Consider also the obvious projections
$$p:\widetilde F_\Vb\to F_\Vb,\quad q:Z_\Vb\to F_\Vb\times F_\Vb.$$
Note that (1.3), (1.5) yield
$$\Fb_\yb=\Ext^*_{G_\Vb}(\kb,\Lc_\yb)=H^*_{G_\Vb}(E_\Vb,\Lc_\yb)=
H^*_{G_\Vb}(\widetilde F_\yb,\kb).$$ Recall the following
isomorphism, see (1.4)
$$H^*_{G_\Vb}(\widetilde F_\yb,\k)=
H^{*}_{G_\Vb}(\widetilde F_\yb,\Dc)[-2d_\yb]=
H_{*}^{G_\Vb}(\widetilde F_\yb,\k)[-2d_\yb].\leqno(1.8)$$ We
obtain a graded $\Sb_\Vb$-module isomorphism
$$\Fb_\yb=
H^{G_\Vb}_*(\widetilde F_\yb,\kb)[-2d_\yb].\leqno(1.9)$$
Next, we equip the $\Sb_\Vb$-module $H_*^{G_\Vb}(Z_\Vb,\k)$ with the
convolution product $\star$ relative to the closed embedding
$Z_\Vb\subset\widetilde F_\Vb\times\widetilde F_\Vb$. See \cite{CG,
sec.~8.6} for details. We obtain a
$\Sb_\Vb$-algebra with 1 which acts on the $\Sb_\Vb$-module
$H^{G_\Vb}_*(\widetilde F_\Vb,\k)$ from the left. The unit of
$H_*^{G_\Vb}(Z_\Vb,\k)$ is the fundamental class of the closed
subvariety $Z^e_\Vb\subset Z_\Vb$. See Section 2.8 below for the
notation.

\proclaim{1.8.~Lemma} (a) The left $H_*^{G_\Vb}(Z_\Vb,\k)$-module
$H^{G_\Vb}_*(\widetilde F_\Vb,\k)$ is faithful.

(b) There is a canonical $\Sb_\Vb$-algebra isomorphism $\Zb_\Vb=
H_*^{G_\Vb}(Z_\Vb,\k)$ such that (1.9) intertwines the
$\Zb_\Vb$-action on $\Fb_\Vb$ and the $H^{G_\Vb}_*(Z_\Vb,\k)$-action
on $H^{G_\Vb}_*(\widetilde F_\Vb,\k)$.

\endproclaim

\noindent{\sl Proof :} Part $(b)$ is proved as in \cite{CG,
Prop.~8.6.35}. Note that in loc.~cit.~both claims are proved for
non-equivariant homology. However the proof uses only standard tools
and generalizes to the equivariant setting. Part $(a)$ is standard.
A proof of $(a)$ is given in Remark 2.21 below for the comfort of
the reader.

\qed

\vskip3mm

\subhead 1.9.~Shift of the grading\endsubhead 
For each $\yb,\yb'\in Y_\nu$ we set
$$\gathered
{}^\delta\Zb_{\yb,\yb'}=
\Ext^*_{G_\Vb}({}^\delta\!\Lc_\yb,{}^\delta\!\Lc_{\yb'}),
\quad{}^\delta\Fb_\yb=\Ext^*_{G_\Vb}({}^\delta\!\Lc_\yb,\Dc).
\endgathered
\leqno(1.10)$$ Consider the graded $\Sb_\Vb$-modules
$$\gathered
{}^\delta\Zb_\Vb=
\Ext^*_{G_\Vb}({}^\delta\!\Lc_\Vb,{}^\delta\!\Lc_{\Vb}),
\quad{}^\delta\Fb_\Vb=\Ext^*_{G_\Vb}({}^\delta\!\Lc_\Vb,\Dc), \cr
\Lc_\Vb=\bigoplus_{\ib\in I^\nu}\Lc_{\ib},\quad
{}^\delta\!\Lc_\Vb=\bigoplus_{\ib\in I^\nu}{}^\delta\!\Lc_{\ib}.
\endgathered$$
Thus we have
$$
\gathered {}^\delta\Zb_\Vb=\bigoplus_{\ib,\ib'\in
I^\nu}{}^\delta\Zb_{\ib,\ib'},
\quad{}^\delta\Fb_\Vb=\bigoplus_{\ib\in I^\nu}{}^\delta\Fb_\ib, \cr
{}^\delta\Zb_{\ib,\ib'}= \Zb_{\ib,\ib'}[d_{\ib}-d_{\ib'}], \quad
{}^\delta\Fb_{\ib}= \Fb_{\ib}[d_\ib].
\endgathered\leqno(1.11)$$

\vskip3mm

\subhead 1.10.~Notation\endsubhead The symbol $\star$ will denote both the multiplication in the
algebras $\Zb_\Vb$, $H_*^{G_\Vb}(Z_\Vb,\k)$ and the left action on
 $\Fb_\Vb$, $H^{G_\Vb}_*(\widetilde
F_\Vb,\k)$. Let 1 denote the unit in $\Zb_\Vb$ and
$H_*^{G_\Vb}(Z_\Vb,\k)$. In this paper we'll never consider the homological grading
on $H_*^{G_\Vb}(Z_\Vb,\k)$ because the isomorphism $\Zb_\Vb=H_*^{G_\Vb}(Z_\Vb,\k)$
in Lemma 1.8$(b)$ is not homogeneous of degree 0. The gradings on $\Zb_\Vb$ and ${}^\delta\Zb_\Vb$
will always be the ones given in Sections 1.3 and 1.9.

\subhead 2. Computation of the algebra $\Zb_\Vb$\endsubhead

This section contains background material needed to compute the
graded Ext-algebras $\Zb_\Vb$, ${}^\delta\Zb_\Vb$ 
introduced in the first section. More precisely we identify
$\Zb_\Vb$, $\Fb_\Vb$ with $H_*^{G_\Vb}(Z_\Vb,\kb)$, 
$H_*^{G_\Vb}(\tilde F_\Vb,\k)$ via Lemma 1.8 and (1.9) and we compute 
explicitely the $H_*^{G_\Vb}(Z_\Vb,\kb)$-action on
$H_*^{G_\Vb}(\tilde F_\Vb,\k)$. We also consider the grading in Corollary 2.25.

\subhead{2.1.~Notation}\endsubhead In this section we fix once for all a
nonzero element $\nu\in\NN I$ and an object $\Vb\in\Vcb_\nu$. Fix
also a maximal torus $T_\Vb\subset G_\Vb$. We'll abbreviate
$G=\GL(\Vb)$ and $F=F(\Vb)$. It is convenient to see $F_\Vb$ as the
closed subvariety of $F$ consisting of all flags which are fixed
under the action of the center of $G_\Vb$.

Let $\Sen_\Vb$, $\Sen$ be the Weyl groups of the pairs
$(G_\Vb,T_\Vb)$, $(G,T_\Vb)$. There is a canonical embedding
$\Sen_\Vb\subset\Sen$. The group $\Sen$ acts freely transitively on
the set $F_\Vb^{T_\Vb}$ of the flags which are fixed by the
$T_\Vb$-action. We'll write $e$ for the unit in $\Sen$.

\subhead{2.2.~The $\Sen$-action on $T_\Vb$-fixed flags}\endsubhead 
Fix once for all a flag $\phi_\Vb\in
F_\Vb^{T_\Vb}$. Write $\phi_{\Vb,w}=w(\phi_\Vb)$ for each
$w\in\Sen$. Let $\ib_w$ be the unique sequence such that
$\phi_{\Vb,w}\in F_{\ib_w}$. There is bijection
$$\Sen_\Vb\setminus\Sen\to I^\nu,\quad\Sen_\Vb w\mapsto
\ib_w.$$ For each  $\ib\in I^\nu$ we have $$\widetilde
F_\ib^{T_\Vb}\simeq F_\ib^{T_\Vb}=\{\phi_{\Vb,w};w\in\Sen_\ib\},$$
where $\Sen_\ib$ is the right $\Sen_\Vb$-coset
$$\Sen_\ib=\{w\in\Sen;\ib_w=\ib\}.$$

Let $B_{\Vb,w}$ be the stabilizer of the flag $\phi_{\Vb,w}$ under
the $G_\Vb$-action. It is a Borel subgroup of $G_\Vb$ containing
$T_\Vb$. Let $N_{\Vb,w}$ be the unipotent radical of $B_{\Vb,w}$.
We'll abbreviate $$F_w=F_{\ib_w},\ \Sen_w=\Sen_{\ib_w},\
\pi_w=\pi_{\ib_w},\ d_w=d_{\ib_w},\ B_\Vb=B_{\Vb,e},\
N_\Vb=N_{\Vb,e}.$$ Note that $\Sen_w=\Sen_\Vb w$ and that we have an
isomorphism of $G_\Vb$-varieties
$$G_\Vb/B_{\Vb,w}\to F_{w},\quad g\mapsto g\phi_{\Vb,w}.$$

\subhead{2.3.~Identification of $\Sen$ with the symmetric group}\endsubhead 
Set $m=|\nu|$, a positive integer.
Recall that $\ib_e$ is the unique sequence such that the flag $\phi_\Vb$
belongs to $F_{\ib_e}$. Write $\ib_e=(i_1,i_2,\dots i_m)$. 
Fix one-dimensional $T_\Vb$-submodules
$\Db_1,\Db_2,\dots \Db_m\subset\Vb$ such that
$$\phi_\Vb=(\Vb\supset \Db_1\oplus \Db_2\oplus\cdots
\Db_{m-1}\supset\cdots\supset \Db_1\oplus \Db_2 \supset \Db_1\supset
0).$$  Then $\Db_l\subset\Vb_{i_l}$ for each $l$. There is a
canonical group isomorphism $\Sen=\Sen_m$ such that
$w(\Db_l)=\Db_{w(l)}$ for each $w,l$. Let 
$\Sen_\nu$ be the image in $\Sen_m$ of the subgroup
$\Sen_\Vb$ of $\Sen$.
The symmetric group $\Sen_m$ acts on the set $I^\nu$ by permutations
: view a sequence $\ib$ as the map $\{1,2,\dots m\}\to I$, $l\mapsto
i_l$ and set $w(\ib)=\ib\circ w^{-1}$ for $w\in\Sen_m$. Then we have
$$\Sen_\nu=\{w\in\Sen_m;w(\ib_e)=\ib_e\}.$$
Under the canonical isomorphism $\Sen=\Sen_m$ we have
$w(\ib_e)=\ib_{w^{-1}}$ for each $w$.

\subhead{2.4.~The root system}\endsubhead Let $B$ be the stabilizer of the flag
$\phi_\Vb$ in $G$. Note that $B_\Vb=B\cap G_\Vb$. Let $\Delta$ be
the set of roots of $(G,T_\Vb)$ and $\Delta^+\subset\Delta$ be the
set of positive roots relative to the Borel subgroup $B\subset G$.
We abbreviate $\Delta^-=-\Delta^+$. Let $\Pi\subset\Sen$ be  the set
of simple reflections. Let $\leqslant$ and $w\mapsto\ell(w)$ denote
the Bruhat order and the length function on $\Sen_m$ or $\Sen$. The
Weyl group $\Sen$ acts on $\Delta$.
Let $\Delta_\Vb$ be the set of roots of $(G_\Vb,T_\Vb)$. Write
$\Delta_\Vb^+=\Delta^+\cap\Delta_\Vb$ and
$\Delta_\Vb^-=-\Delta_\Vb^+$. The action of $\Sen_\Vb$ on $\Delta$
preserves the subset $\Delta_\Vb$. Let
$\chi_1,\chi_2,\dots\chi_m\in\ten_\Vb^*$ be the weights of the lines
$\Db_1,\Db_2,\dots \Db_m$ respectively. The set of simple roots in $\Delta^+$ is
$$\{\chi_l-\chi_{l+1};\,l=1,2,\dots,m-1\}.$$
Let $s_l\in\Pi$ denote the reflection with respect to the simple root $\chi_l-\chi_{l+1}$.
Under the identification $\Sen=\Sen_m$ the reflection $s_l$ is taken to the simple transposition $(l,l+1)$.

\subhead{2.5.~The stratification of $F_\Vb\times F_\Vb$}\endsubhead 
The group $G$ acts diagonally on $F\times
F$. The action of the subgroup $G_\Vb$ preserves the subset
$F_\Vb\times F_\Vb$. For $x\in\Sen$ let $O^{x}_\Vb$ be the set of
all pairs of flags in $F_\Vb\times F_\Vb$ which are in relative
position $x$. More precisely, write
$$\phi_{\Vb,w',w}=(\phi_{\Vb,w'},\phi_{\Vb,w}),\quad\forall
w,w'\in\Sen.$$ Then we set $$O^{x}_\Vb=(F_\Vb\times F_\Vb)\cap
G\phi_{\Vb,e,x}.$$ Let $\bar O_\Vb^{x}$ be the Zariski closure of
$O_\Vb^{x}$. It is the set of pairs of flags in $F_\Vb\times F_\Vb$
which are in relative position $\leqslant x$. For any $w',w\in\Sen$
we write also $$O^{x}_{w',w}= O^{x}_\Vb\cap( F_{w'}\times
F_{w}),\quad\bar O^{x}_{w',w}=\bar O^{x}_\Vb\cap( F_{w'}\times
F_{w}).$$
Fix a simple reflection $s\in\Pi$. Put
$$P_{\Vb,w,ws}=B_{\Vb,w}\{wsw^{-1},1\}B_{\Vb,w}.$$
It is a parabolic subgroup of $G_\Vb$ containing $B_{\Vb,w}$. The
following lemma is standard. Its proof is left to the reader.

\proclaim{2.6.~Lemma} (a) The set of elements of $O^x_\Vb$ which are
fixed under the diagonal action of $T_\Vb$ is
$\{\phi_{\Vb,w,wx};w\in\Sen\}$.

(b) The variety $\bar O_\Vb^{s}$ is smooth and is equal to
$O_\Vb^{s}\cup O_\Vb^{e}$.

(c) We have $\bar O^{s}_{w,w'}=\emptyset$ unless $w'=w$ or $ws$.

(d) Assume that $ws\notin\Sen_\Vb w$. We have $$F_{ws}\neq F_w,\quad
B_{\Vb,ws}=B_{\Vb,w},\quad\bar O^s_{w,ws}=O^s_{w,ws},\quad\bar
O^s_{w,w}=O^e_{w,w}.$$

(e) Assume that $ws\in\Sen_\Vb w$. We have
$$O^s_{w,ws}=O^s_{w,w},\quad F_{ws}=F_w,\quad B_{\Vb,ws}\neq
B_{\Vb,w}.$$ There is an isomorphism of $G_\Vb$-varieties
$$G_\Vb\times_{B_{\Vb,w}}(P_{\Vb,w,ws}/B_{\Vb,w})\to\bar
O^s_{w,w},\quad (g,h)\mapsto(g\phi_{\Vb,w},gh\phi_{\Vb,w}).$$

\endproclaim

\subhead 2.7.~Example\endsubhead Set $I=\{i,j\}$, $\nu=i+2\,j$,
$\Vb_i=\Db_1$, $\Vb_j=\Db_2\oplus \Db_3$. So
$I^\nu=\{\ib_e,\ib_{s},\ib_{t}\}$ where $s=s_1$, $t=s_1s_2s_1$, and $\ib_e=(i,j,j)$, $\ib_{s}=(j,i,j)$,
$\ib_{t}=(j,j,i)$. We have $$\aligned &\phi_\Vb=(\Vb\supset
\Db_1\oplus \Db_2\supset \Db_1\supset 0), \hfill\cr
&\phi_{\Vb,s}=(\Vb\supset \Db_1\oplus \Db_2\supset \Db_2\supset
0),\hfill\cr &\phi_{\Vb,t}=(\Vb\supset \Db_2\oplus \Db_3\supset
\Db_3\supset 0).\endaligned$$ We have also
$$\aligned
&F_e=\{\Vb\supset \Db_1\oplus \Db\supset \Db_1\supset 0; \,
\Db\in\PP(\Vb_j)\},\hfill\cr &F_{s}=\{\Vb\supset \Db_1\oplus
\Db\supset \Db\supset 0; \, \Db\in\PP(\Vb_j)\}, \hfill\cr
&F_{t}=\{\Vb\supset\Vb_j\supset\Db\supset 0; \, \Db\in\PP(\Vb_j)\}.
\hfill\cr
\endaligned$$
Finally we have $$\bar O^{s}_{s,t}=\bar O^{s}_{1,t}=\emptyset,\quad
\bar O^{s}_{t,t}=F_t\times F_t, \quad\bar
O^{s}_{s,s}=O^{e}_{s,s}=\Delta_{F_{s}}, \quad\bar
O^{s}_{e,e}=O^{e}_{e,e}=\Delta_{F_e},$$ where $\Delta$ is the
diagonal, and $$\bar
O^{s}_{e,s}=O^{s}_{e,s}=\{(\Vb\supset\Db_1\oplus\Db\supset\Db\supset
0,\, \Vb\supset\Db_1\oplus\Db\supset\Db_1\supset
0);\Db\in\PP(\Vb_j)\}.$$

\subhead 2.8.~The stratification of $Z_\Vb$ and $\Zb_\Vb$\endsubhead 
Recall the notation in Section 1.7. For
$x\in\Sen$ let $Z_\Vb^{x}\subset Z_\Vb$ be the Zariski closure of
the locally closed subset $q^{-1}(O^{x}_\Vb)$. Put
$$Z_\Vb^{\leqslant x}=\bigcup_{w\leqslant x}Z_\Vb^w,\quad Z_{\ib',\ib}^{\leqslant x}=Z_\Vb^{\leqslant
x}\cap Z_{\ib',\ib},\quad\forall\,\ib,\ib'\in I^\nu.$$ Hence
$Z_\Vb^x$, $Z_\Vb^{\leqslant x}$, $Z_{\ib',\ib}^{\leqslant x}$ are
closed $G_\Vb$-subvarieties of $Z_\Vb$. Lemma 1.8$(b)$ yields
$\kb$-vector space isomorphisms
$$\Zb_{\ib',\ib}=H_*^{G_\Vb}(Z_{\ib',\ib},\k),\quad\forall\,\ib,\ib'\in I^\nu.$$
We set $$\Zb_\Vb^{\leqslant x}=H_*^{G_\Vb}(Z_\Vb^{\leqslant
x},\k),\quad \Zb_\Vb^{e}=\Zb_\Vb^{\leqslant e},
\quad\Zb_{\ib',\ib}^{\leqslant w}= \Zb_\Vb^{\leqslant
w}\cap\Zb_{\ib',\ib}.$$

\proclaim{2.9.~Lemma} (a) The direct image by the closed embedding
$Z_\Vb^{\leqslant x}\subset Z_\Vb$ gives an injective left graded
$\Sb_\Vb$-module homomorphism $\Zb_\Vb^{\leqslant x}\to\Zb_\Vb$.

(b) For each $x,y\in\Sen$ such that $\ell(xy)=\ell(x)+\ell(y)$ the
convolution product in $\Zb_\Vb$ yields an inclusion
$\Zb_\Vb^{\leqslant x}\star\Zb_\Vb^{\leqslant
y}\subset\Zb_\Vb^{\leqslant xy}$.

(c) We have $1_\ib\in\Zb_{\ib,\ib}^{e}$ for each  $\ib\in I^\nu$.

\endproclaim

\noindent{\sl Proof :} Parts $(a)$, $(b)$ are standard, see
\cite{CG, chap.~6} for the case of the equivariant K-theory of the
Steinberg variety. The proof is the same in our case. Part $(c)$ is
obvious : since the convolution product by $1_\ib$ is the identity
$\Zb_{\ib,\ib}\to\Zb_{\ib,\ib}$ we must have
$1_{\ib}=[Z^e_{\ib,\ib}]$.

\qed

\subhead{2.10.~Euler classes in $\Pb_\Vb$}\endsubhead 
Consider the $\kb$-algebra 
$$\Pb_\Vb=\Sb_{T_\Vb}\leqno(2.1)$$
Assume that $M$ is a finite dimensional
representation of the Lie algebra $\ten_\Vb$. 
For each linear form $\l\in\ten_\Vb^*$ let $M[\l]\subset M$ be the weight
subspace associated with $\l$. Let $\ch(M)=\sum_\l\dim(M[\l])\l$
be the character of $M$. Let $\Eu(M)$ be the
determinant of $M$, viewed as a homogeneous polynomial of degree
$\dim(M)$ on $\ten_\Vb$, i.e., an element of degree $2\dim(M)$ in
$\Pb_\Vb$. If $M$ is a finite dimensional representation of $T_\Vb$
we write $\Eu(M)$ again for the polynomial associated to the
differential of $M$, a module over $\ten_\Vb$.
Now, assume that $X$ is a quasi-projective $T_\Vb$-variety and that
$x\in X^{T_\Vb}$ is a smooth point of $X$. The cotangent space
$T^*_xX$ at $x$ is equipped with a natural representation of
$T_\Vb$. We'll abbreviate $\Eu(X,x)=\Eu(T_x^*X)$. 
We'll be particularly interested by the following element
$$\Lambda_{w}=\Eu(\widetilde F_\Vb,\phi_{\Vb,w}),\quad\forall w\in \Sen.$$
Note that $\Lambda_w$ is an element of $\Pb_\Vb$ of degree $2d_w$.

\subhead 2.11.~Notation\endsubhead For each $w\in\Sen$ let
$$\een_{\Vb,w}=\{x\in E_\Vb;\phi_{\Vb,w}\ \text{is}\
x\text{-stable}\}.$$ Under restriction the natural $G_\Vb$-action on
$E_\Vb$ yields a representation of $B_{\Vb,w}$ on $\een_{\Vb,w}$.
There is an isomorphism of $G_\Vb$-varieties
$$G_\Vb\times_{B_{\Vb,w}}\een_{\Vb,w}\to\widetilde F_{w},
\quad(g,x)\mapsto (g\phi_{\Vb,w},gx).$$ Under this isomorphism the
map $\pi_{w}:\widetilde F_{w}\to E_\Vb$ is identified with the map
$$G_\Vb\times_{B_{\Vb,w}}\een_{\Vb,w}\to E_\Vb,\quad
(g,x)\mapsto gx.$$ As a $T_\Vb$-module $\een_{\Vb,w}$ is the sum of
the weight subspaces of $E_\Vb$ whose weights belong to
$w(\Delta^+)$. For $w,w'\in\Sen$ we write
$$\een_{\Vb,w,w'}=\een_{\Vb,w}\cap\een_{\Vb,w'},\quad
\den_{\Vb,w,w'}=\een_{\Vb,w}/\een_{\Vb,w,w'}.$$

\subhead{2.12}\endsubhead Fix a simple reflection $s\in\Pi$. Let $\gen_\Vb$, $\ten_\Vb$,
$\nen_{\Vb,w}$, $\pen_{\Vb,w,ws}$ be the Lie algebras of $G_\Vb$,
$T_\Vb$, $N_{\Vb,w}$, $P_{\Vb,w,ws}$ respectively. Let
$\nen_{\Vb,w,ws}$ be the nilpotent radical of $\pen_{\Vb,w,ws}$. We
have
$$\nen_{\Vb,w,ws}=\nen_{\Vb,w}\cap\nen_{\Vb,ws}.$$ Let
$\men_{\Vb,w,ws}=\nen_{\Vb,w}/\nen_{\Vb,w,ws}$. So we have the
following $T_\Vb$-module isomorphisms
$$\nen_{\Vb,w}=\nen_{\Vb,w,ws}\oplus\men_{\Vb,w,ws},\quad
\men_{\Vb,w,ws}\simeq\men_{\Vb,ws,w}^*.$$  Note that
$\nen_{\Vb,w}\subset\gen_\Vb$ is the sum of the root subspaces whose
weights belong to $w(\Delta^+)\cap\Delta_\Vb$.

\subhead{2.13.~Reduction to the torus}\endsubhead Recall the graded $\k$-algebra
$\Pb_\Vb=\Sb_{T_\Vb}$. The canonical action of
$\Sen$ on $T_\Vb$ yields a $\Sen$-action on $\Pb_\Vb$. The
restriction of functions from $\gen_\Vb$ to $\ten_\Vb$ gives an
isomorphism of graded $\k$-algebras
$$\Sb_\Vb=(\Pb_\Vb)^{\Sen_\Vb}.\leqno(2.2)$$
Now, recall the notation in Sections 2.3-4. So
$\chi_1,\chi_2,\dots\chi_m$ are the weights of the lines
$\Db_1,\Db_2,\dots \Db_m$ respectively. These weights generate the
algebra $\Pb_\Vb$ and they have the degree 2. Under the
identification $\Sen\simeq\Sen_m$ the action of $w$ on $\Pb_\Vb$ is
given by
$$f=f(\chi_1,\cdots,\chi_m)\mapsto w(f)=f(\chi_{w(1)},\dots,\chi_{w(m)}).$$
Next, we recall here a standard result for a future use. If $X$ is a
quasi-projective $G_\Vb$-variety then the $\Pb_\Vb$-module
$H_*^{T_\Vb}(X,\k)$ is equipped with a canonical
$\Pb_\Vb$-skewlinear action of the group $\Sen_\Vb$. It is
well-known that the forgetful map gives a $\Sb_\Vb$-module
isomorphism
$$H_*^{G_\Vb}(X,\k)\to H_*^{T_\Vb}(X,\k)^{\Sen_\Vb}.\leqno(2.3)$$

\subhead{2.14.~The $\kb$-algebra structure on $\Fb_\Vb$}\endsubhead 
Let $\ib\in I^\nu$. Assigning to a formal
variable $x_\ib(l)$ of degree 2 the first equivariant Chern class of
the $G_\Vb$-equivariant line bundle $\Oc_{\widetilde F_\ib}(l)$ in
Section 1.1 for each $l=1,2,\dots m$ yields a graded $\k$-algebra
isomorphism
$$\k[x_\ib(1),x_\ib(2),\dots x_\ib(m)]=
H^*_{G_\Vb}(\widetilde F_\ib,\k).$$ Composing this map by the
isomorphisms (1.8), (1.9) we get a canonical isomorphism of graded
$\k$-vector spaces
$$\k[x_\ib(1),x_\ib(2),\dots x_\ib(m)]
=H_*^{G_\Vb}(\widetilde F_\ib,\k)[-2d_\ib]=\Fb_\ib.\leqno(2.4)$$ The
multiplication of polynomials equip each $\Fb_\ib$ with an obvious
structure of graded $\k$-algebra. Thus $\Fb_\Vb$ is also a graded
$\k$-algebra by (1.7).

\subhead{2.15.~The $\Sen_m$-action and the $\Sb_\Vb$-action on $\Fb_\Vb$}
\endsubhead 
For each $\ib\in I^\nu$ and each
$w\in\Sen_\ib$ the pull-back by the closed embedding
$\{\phi_{\Vb,w}\}\subset \widetilde F_\ib$ yields a graded
$\k$-algebra isomorphism
$$\Fb_\ib\to\Pb_\Vb,\quad
f(x_\ib(1),\dots,x_\ib(m))\mapsto f(\chi_{w(1)},\dots\chi_{w(m)}).\leqno(2.5)$$ 
We'll write $w(f)$ for the right hand side.
This isomorphism is
not canonical, because it depends on the choice of $w$. 

Consider the $\Sen_m$-action on $\Fb_\Vb$ given by 
$$w\Fb_\ib=\Fb_{w(\ib)},\quad w f(x_\ib(1),\dots
x_\ib(m))=f(x_{w(\ib)}(w(1)),\dots,x_{w(\ib)}(w(m))).$$ 

Next, consider the canonical $\Sb_\Vb$-action on $\Fb_\Vb$ coming from the 
$G_\Vb$-equivariant cohomology. It can be regarded as a $\Sb_\Vb$-action on
$\bigoplus_\ib\k[x_\ib(1),x_\ib(2),\dots x_\ib(m)]$ which is described in the 
following way.
The composition of the obvious projection $\Fb_\Vb\to\Fb_\ib$ with the map (2.5)
identifies the graded $\k$-algebra of the $\Sen_m$-invariant polynomials in the
$x_\ib(l)$'s with $(\Pb_\Vb)^{\Sen_\Vb}=\Sb_\Vb$.
This isomorphism is canonical.
The $\Sb_\Vb$-action on $\Fb_\Vb$ is the composition of this isomorphism and of 
the multiplication by $\Sen_m$-invariant polynomials.

\subhead{2.16.~Examples}\endsubhead $(a)$ Fix a simple reflection
$s_l\in\Pi$. For each $w\in\Sen$ the $T_\Vb$-module
$\nen_{\Vb,w}$ is the sum of the root subspaces of weight
$\chi_{w(k)}-\chi_{w(k')}$ with $k<k'$ which are contained in
$\gen_\Vb$. The following hold

\item{$\bullet$} if $ws_l\in\Sen_\Vb w$ then we have
$$\Eu(\nen_{\Vb,ws_l})=-\Eu(\nen_{\Vb,w}),\quad
\Eu(\men_{\Vb,w,ws_l})=-\Eu(\men_{\Vb,ws_l,w})=\chi_{w(l)}-\chi_{w(l+1)},$$

\item{$\bullet$}
if $ws_l\notin\Sen_\Vb w$ then we have
$$\nen_{\Vb,ws_l}=\nen_{\Vb,w},\quad
\Eu(\men_{\Vb,w,ws_l})=\Eu(\men_{\Vb,ws_l,w})=0.$$

$(b)$ Let $s_l$ be as above. Given $w\in\Sen$ we write
$$\ib_w=(i_1,i_2,\dots i_m),
\quad \Vb^k_w=\Db_{w(1)}\oplus\Db_{w(2)}\oplus\cdots
\Db_{w(k)},\quad\forall k.$$ Note that in Section 2.3 we used the
(different) notation $\ib_e=(i_1,i_2,\dots i_m)$. The following
properties are straightforward
$$\aligned
&\phi_{\Vb,w}=(\Vb=\Vb_w^m\supset\Vb^{m-1}_w\supset\cdots\supset\Vb_w^0=0),\hfill\cr
&\Db_{w(k)}\subset \Vb_{i_k},\hfill\cr &\Vb^k_w=\Vb_{ws_l}^k\
\roman{if}\ k\neq l,\hfill\cr &x(\Vb^k_w)\subset
\Vb^{k-1}_w,\quad\forall x\in \een_{\Vb,w}.
\endaligned\leqno(2.6)$$
The last claim follows from the inclusion $x(\Vb^k_w)\subset
\Vb^{k}_w$, because else $x$ would give a non-zero map
$\Db_{w(k)}\to\Db_{w(k)}$ and this is impossible because
$i_{k}\not\to i_{k}$. So we have
$$\aligned
&\een_{\Vb,w}=\{x\in E_\Vb;\,x(\Vb_w^k)\subset\Vb_w^{k-1},\,\forall
k\}, \cr &\een_{\Vb,w,ws}=\{x\in\een_{\Vb,w};\,
x(\Db_{w(l+1)})\subset\Vb_w^{l-1}\}.
\endaligned$$
Therefore, composing the map
$$\een_{\Vb,w}\to\bigoplus_{h:i_{l+1}\to i_l}\Hom(\Db_{w(l+1)},\Vb_w^l),\quad
x\mapsto (x_h|_{\Db_{w(l+1)}};h\in H_{i_{l+1},i_l}),
$$
with the isomorphism $\Vb_w^l/\Vb_w^{l-1}=\Db_{w(l)}$ we get a
$T_\Vb$-module isomorphism
$$\den_{\Vb,w,ws}=\bigoplus_{h:i_{l+1}\to i_l}\Db_{w(l+1)}^*\otimes\Db_{w(l)}.
\leqno(2.7)$$ Hence we have
$$\eu(\den_{\Vb,w,ws})=(\chi_{w(l)}-\chi_{w(l+1)})^{h_{i_{l+1},i_l}}.
\leqno(2.8)$$

\subhead{2.17.~Localization to the $T_\Vb$-fixed points}\endsubhead 
For each $w\in\Sen$ we set
$\psi_{w}=[\{\phi_{\Vb,w}\}]$, a class in $H_*^{T_\Vb}(\widetilde
F_\Vb,\k)$. For each $w,w'\in\Sen$ we set also
$\psi_{w,w'}=[\{\phi_{\Vb,w,w'}\}]$, a class in
$H_*^{T_\Vb}(Z_{\Vb},\k)$. Let $\Kb_\Vb$ be the fraction field of
$\Pb_\Vb$. For each quasi-projective $T_\Vb$-variety $X$ we'll
abbreviate
$$\Hc_*(X,\k)=
H_*^{T_\Vb}(X,\k)\otimes_{\Pb_\Vb}\Kb_\Vb.$$ We'll also abbreviate
$\psi_w$, $\psi_{w,w'}$ for the elements
$$\psi_w\otimes 1\in\Hc_*(\widetilde F_\Vb,\k),\quad
\psi_{w,w'}\otimes 1\in\Hc_*(Z_{\Vb},\k).$$ The following lemma is
standard. Its proof is left to the reader.

\proclaim{2.18. Lemma} (a) The $\Pb_\Vb$-modules
$H_*^{T_\Vb}(\widetilde F_\Vb,\k)$, $H_*^{T_\Vb}(Z_\Vb,\k)$ are
free.

(b) We have $\Kb_\Vb$-vector spaces isomorphisms
$$\Hc_*(\widetilde
F_\Vb,\k)=\bigoplus_{x\in\Sen}\Kb_\Vb\psi_{x}, \quad
\Hc_*(Z_{\Vb},\k)= \bigoplus_{x,y\in\Sen}\Kb_\Vb\psi_{x,y}.$$

(c) The canonical $\Sen_\Vb$-action on the $\Pb_\Vb$-modules
$H_*^{T_\Vb}(\widetilde F_\Vb,\k)$, $H_*^{T_\Vb}(Z_\Vb,\k)$ (see Section $2.13$) is given
by $w(\psi_{x})=\psi_{wx}$ and $w(\psi_{x,y})=\psi_{wx,wy}$ for each
$w\in\Sen_\Vb$, $x,y\in\Sen$.

(d) Fix a sequence $\ib\in I^\nu$. Composing (2.3), (2.4) and the
obvious map $$H_*^{T_\Vb}(\widetilde F_\ib,\k)\to\Hc_*(\widetilde
F_\ib,\k)$$ we get a canonical embedding
$$\k[x_\ib(1),\dots x_\ib(m)]\to\bigoplus_{w\in\Sen_\ib}\Kb_\Vb\psi_{w},\quad
f(x_\ib(1),\dots
x_\ib(m))\mapsto\sum_{w\in\Sen_\ib}w(f)\Lambda_w^{-1}\psi_w.$$

\endproclaim

\subhead 2.19\endsubhead Fix a simple reflection $s\in\Pi$. Recall the closed $G_\Vb$-subvariety
$Z_\Vb^{s}\subset Z_\Vb$ from Section 2.8. There are unique rational
functions $\Lambda_{w,w'}^{s}\in\Kb_\Vb$, $w,w'\in\Sen$, such that
the image of $[Z_\Vb^s]$ in $\Hc_*(Z_\Vb,\k)$ is
$$[Z_\Vb^{s}]=\sum_{w,w'}\Lambda_{w,w'}^{s}\,\psi_{w,w'}.\leqno(2.9)$$
Note that $\Lambda^s_{w,w'}=0$ unless $w'=w$ or $ws$ by Lemma
2.6$(c)$.
Let $\star$ denote also the convolution products
$$\aligned
&H_*^{T_\Vb}(Z_{\Vb},\k)\times H_*^{T_\Vb}(Z_{\Vb},\k)\to
H_*^{T_\Vb}(Z_{\Vb},\k),\hfill\cr &H_*^{T_\Vb}(Z_{\Vb},\k)\times
H_*^{T_\Vb}(\widetilde F_{\Vb},\k)\to H_*^{T_\Vb}(\widetilde
F_{\Vb},\k)\endaligned$$ relative to the closed embedding
$Z_\Vb\subset\widetilde F_\Vb\times\widetilde F_\Vb$.
Compare Section 1.7.

\proclaim{2.20. Lemma} (a) The forgetful maps below commute with
$\star$
$$\Zb_{\Vb}\to H_*^{T_\Vb}(Z_{\Vb},\k),\quad
\Fb_{\Vb}\to H_*^{T_\Vb}(\widetilde F_{\Vb},\k).$$

(b) For $w,w',w''\in\Sen$ we have
$$\psi_{w',w}\star\psi_{w}=\Lambda_{w}\,\psi_{w'},
\quad\psi_{w'',w'}\star\psi_{w',w}=\Lambda_{w'}\,\psi_{w'',w}.$$

(c) For $w\in\Sen$ we have

$$\aligned
&\ \Lambda_w={\Eu}(\een_{\Vb,w}^*\oplus\nen_{\Vb,w}), \hfill\cr
&\left.\aligned &\Lambda^{s}_{w,w}=
\Eu(\een_{\Vb,ws,w}^*\oplus\nen_{\Vb,w}\oplus\men_{\Vb,w,ws})^{-1}
\hfill\cr &\Lambda^{s}_{w,ws}=
\Eu(\een_{\Vb,ws,w}^*\oplus\nen_{\Vb,w}\oplus\men_{\Vb,ws,w})^{-1}
\endaligned
\right\} \quad\roman{if}\ ws\in\Sen_\Vb w, \hfill\cr &\
\Lambda^s_{w,ws}=
\Lambda^{s}_{w,w}=\Eu(\een_{\Vb,ws,w}^*\oplus\nen_{\Vb,w})^{-1}
\hskip11mm\roman{if}\ ws\notin\Sen_\Vb w.
\endaligned$$

\endproclaim

\noindent{\sl Proof :} Parts $(a)$, $(b)$ are well-known. Let us
concentrate on $(c)$. The fiber at $\phi_{\Vb,w}$ of the vector
bundle $p:\widetilde F_\Vb\to F_\Vb$ is isomorphic to $\een_{\Vb,w}$
as a $T_{\Vb}$-module. Thus the cotangent space to $\widetilde
F_\Vb$ at the point $\phi_{\Vb,w}$ is isomorphic to
$\een_{\Vb,w}^*\oplus \nen_{\Vb,w}$ as a $T_\Vb$-module. This yields
the first formula. The variety $Z^s_\Vb$ is smooth. A standard
computation yields the following
$$\Lambda^{s}_{w,w'}=\Eu(Z_\Vb^{s},\phi_{\Vb,w,w'})^{-1}.$$
The fiber at $\phi_{\Vb,w,w'}$ of the vector bundle $q:Z_\Vb^s\to
F_\Vb\times F_\Vb$  is isomorphic to $\een_{\Vb,ws,w}$ if $w'=ws,w$
and to 0 else, as a $T_\Vb$-module. Therefore, we have
$$\Lambda^{s}_{w,w'}=\Eu(\een_{\Vb,ws,w}^*)^{-1}\Eu(\bar
O_\Vb^{s},\phi_{\Vb,w,w'})^{-1}.$$ If $ws\in\Sen_\Vb w$ the
cotangent spaces to the variety $\bar O_\Vb^{s}$ at the points
$\phi_{\Vb,w,ws}$, $\phi_{\Vb,w,w}$ are isomorphic to
$\nen_{\Vb,w}\oplus\men_{\Vb,ws,w}$,
$\nen_{\Vb,w}\oplus\men_{\Vb,w,ws}$ respectively as $T_\Vb$-modules,
by Lemma 2.6$(e)$. If $ws\notin\Sen_\Vb w$ they are both isomorphic
to $\nen_{\Vb,w}$ by Lemma 2.6$(d)$. This yields the remaining
formulas.

\qed

\vskip3mm

\subhead 2.21.~Remark\endsubhead We can now prove Lemma 1.8$(a)$. By
Lemma 2.20$(a)$ and (2.3) it is enough to prove that the convolution
product yields a faithful representation of $H_*^{T_\Vb}(Z_\Vb,\k)$
on $H_*^{T_\Vb}(\widetilde F_\Vb,\k)$. By Lemma 2.18$(a)$ it is thus
enough to prove that $\star$ yields indeed a faithful representation
of $\Hc_*(Z_\Vb,\k)$ on $\Hc_*(\widetilde F_\Vb,\k)$. This is
obvious by Lemma 2.18$(b)$, 2.20$(b)$.

\subhead 2.22.~Description of the $\Zb_\Vb$-action on $\Fb_\Vb$\endsubhead 
Fix a simple reflection $s=s_l\in\Pi$. The fundamental class of
$Z^{s}_\Vb$ yields an element $\sigma(l)\in\Zb_\Vb^{\leqslant s}$.
For each sequences $\ib,\ib'\in I^\nu$ we write
$$\sigma_{\ib',\ib}(l)=1_{\ib'}\star\sigma(l)\star 1_\ib
\in\Zb_\Vb^{\leqslant s}.$$
Next, let $k$ be any positive integer $\leqslant m$. The pull-back
of the first equivariant Chern class of the line bundle
$\bigoplus_\ib\Oc_{\widetilde F_\ib}(k)$ by the obvious map
$$Z_\Vb^e\subset\widetilde F_\Vb\times\widetilde F_\Vb
\to\widetilde F_\Vb$$ belongs to $H^*_{G_\Vb}(Z_\Vb^e,\k)$. It
yields an element $\varkappa(k)\in\Zb_\Vb^e$ as in Section 1.7. For
each sequences $\ib,\ib'\in I^\nu$ we write
$$\varkappa_{\ib',\ib}(k)=1_{\ib'}\star \varkappa(k)\star 1_\ib
\in\Zb_\Vb^e.$$
Now, recall that $\Zb_\Vb^{\leqslant s},\Zb_\Vb^e\subset\Zb_\Vb$ by
Lemma 2.9$(a)$ and that $\Zb_\Vb$ acts on $\Fb_\Vb$ by Section 1.6.
The action of $\sigma_{\ib',\ib}(l)$, $\varkappa_{\ib',\ib}(k)$ on
$\Fb_\Vb$ is given by the following formulas.

\proclaim{2.23. Proposition} Let $k,l$ be as above. Fix sequences
$\ib,\ib',\ib''\in I^\nu$ and fix an element $f\in\Fb_\ib$. Write
$\ib=(i_1,\dots,i_m)$.

(a) We have $1_{\ib'}\star f=f$ if $\ib=\ib'$ and 0 else.

(b) We have $\varkappa_{\ib'',\ib'}(k)\star f=0$ unless
$\ib''=\ib'=\ib$ and $\varkappa_{\ib,\ib}(k)\star f=x_\ib(k)f$.

(c) We have $\sigma_{\ib'',\ib'}(l)\star f=0$ unless $\ib'=\ib$ and
$\ib''\in\{s_l(\ib),\ib\}$. More precisely

\item{$\bullet$} if $s_l(\ib)=\ib$ then
$$\sigma_{\ib,\ib}(l)\star f=(f-s_l(f))/(x_\ib(l)-x_\ib(l+1)),$$
\item{$\bullet$} if $s(\ib)\neq\ib$ then
$$\aligned
&\sigma_{s_l(\ib),\ib}(l)\star f=
(x_{s_l(\ib)}(l+1)-x_{s(\ib)}(l))^{h_{i_l,i_{l+1}}} s_l(f),\cr
&\sigma_{\ib,\ib}(l)\star f=
(x_\ib(l+1)-x_\ib(l))^{h_{i_{l+1},i_{l}}} f.\endaligned$$

\endproclaim

\noindent{\sl Proof :} We'll abbreviate $s=s_l$. Parts $(a)$, $(b)$ are standard and are left
to the reader. Let us concentrate on $(c)$. The first claim is
obvious, because $Z^{s}_{w',w}=\emptyset$ unless $w'=w$ or $ws$ and
$\ib_{ws}=s(\ib_w)$. Now, given $\ib'=\ib$ or $s(\ib)$ we must
compute the linear operator
$$\Fb_\ib\to\Fb_{\ib'},\quad f\mapsto\sigma_{\ib',\ib}(l)\star f.\leqno(2.10)$$
By (2.4) we have a $\kb$-vector space isomorphism
$$\Fb_\ib=\kb[x_\ib(1),\dots x_\ib(m)].$$
So Lemma 2.18$(d)$ yields an embedding
$$\Fb_\ib\to \bigoplus_{w\in\Sen_\ib}\Kb_\Vb\psi_w,
\quad f(x_\ib(1),\dots x_\ib(m))\mapsto f_\ib,\quad
f_\ib=\sum_{w\in\Sen_\ib}w(f)\Lambda_w^{-1}\psi_w,\quad
f\in\Pb_\Vb.$$ Under this inclusion the map (2.10) is given by
$$f_\ib\mapsto\sum_{w'\in\Sen_{\ib'}} g_{w'}\psi_{w'},\quad
 g_{w'}=\sum_{w\in\Sen_\ib}w(f)\Lambda^s_{w',w}.$$ See (2.9)
and Lemma 2.20$(b)$. We claim that the rhs is equal to $g_{\ib'}$
for some polynomial $g\in\Pb_\Vb$ that we'll compute explicitly. So
we must find $g$ such that $$
g_{w'}=w'(g)\Lambda_{w'}^{-1},\quad\forall
w'\in\Sen_{\ib'}.\leqno(2.11)$$ In the rest of the proof we'll
assume that the following hold
$$w\in\Sen_\ib,\ w'\in\Sen_{\ib'},\ w'=w\ \roman{or}\ ws,\
\ib'=\ib\ \roman{or}\ s(\ib).$$ As in Example 2.16$(a)$ we'll write
$$\ib=\ib_w=(i_1,i_2,\dots i_m).$$

{$(i)$} Assume that $s(\ib)=\ib$. Then $\ib'=\ib$, $ws\in\Sen_\Vb w$
and we have $$g_{w'}=w'(f)\Lambda^s_{w',w'}+
w's(f)\Lambda^s_{w',w's}.$$ Further there is no arrow in $H$ joining
$i_l$ and $i_{l+1}$. Thus by (2.7) we have
$$\een_{\Vb,w}=\een_{\Vb,ws},$$ and by Lemma 2.16$(a)$ we have
$$\Eu(\nen_{\Vb,ws})=-\Eu(\nen_{\Vb,w}),\quad
\Eu(\men_{\Vb,w,ws})=-\Eu(\men_{\Vb,ws,w})=\chi_{w(l)}-\chi_{w(l+1)}.$$
By Lemma 2.20$(c)$ we have also
$$
\aligned &\Lambda_{w} ={\Eu}(\een_{\Vb,w}^*\oplus\nen_{\Vb,w})
={\Eu}(\een_{\Vb,ws}^*\oplus\nen_{\Vb,ws}) =-\Lambda_{ws}, \hfill\cr
&\Lambda^{s}_{w,w}
=\Eu(\een_{\Vb,ws,w}^*\oplus\nen_{\Vb,w}\oplus\men_{\Vb,w,ws})^{-1}
=\Lambda_w^{-1}\Eu(\men_{\Vb,w,ws})^{-1}, \hfill\cr
&\Lambda^{s}_{w,ws}=
\Eu(\een_{\Vb,ws,w}^*\oplus\nen_{\Vb,w}\oplus\men_{\Vb,ws,w})^{-1}
=\Lambda_{ws}^{-1}\Eu(\men_{\Vb,w,ws})^{-1}.
\endaligned
$$
Therefore we obtain
$$\aligned
g_{w'}&=w'(f)\Eu(\men_{\Vb,w',w's})^{-1}\Lambda_{w'}^{-1}+
w's(f)\Eu(\men_{\Vb,w',w's})^{-1}\Lambda_{w's}^{-1},\hfill\cr
&=(w'(f)\Lambda_{w'}^{-1}+
w's(f)\Lambda_{w's}^{-1})/(\chi_{w'(l)}-\chi_{w'(l+1)}),\hfill\cr
&=w'(g)\Lambda_{w'}^{-1},
\endaligned$$
where $g=(f-s(f))/(\chi_l-\chi_{l+1}).$

{$(ii)$} Assume that $s(\ib)\neq\ib$, i.e., that $ws\notin\Sen_\Vb
w$. Then one of the two following alternative holds :

\vskip1mm

\item{$\bullet$}
$\ib'=s(\ib)$, $w'=ws$ and $g_{w'}=
w's(f)(\Lambda^s_{ws,w}\Lambda_{ws})\Lambda_{w'}^{-1}$,

\vskip1mm

\item{$\bullet$}
$\ib'=\ib$, $w'=w$ and $g_{w'}=
w'(f)(\Lambda^s_{w,w}\Lambda_{w})\Lambda_{w'}^{-1}$.

\vskip1mm

\noindent Now, by (2.8) we have
$$\aligned
&\Eu(\den_{\Vb,ws,w})=(\chi_{w(l+1)}-\chi_{w(l)})^{h_{i_l,i_{l+1}}},\hfill\cr
&\Eu(\den_{\Vb,w,ws})=(\chi_{w(l)}-\chi_{w(l+1)})^{h_{i_{l+1},i_l}}.
\endaligned$$ By Lemma 2.16$(a)$ we
have
$$\Eu(\nen_{\Vb,ws})=\Eu(\nen_{\Vb,w}).$$
By Lemma 2.20$(c)$ we have
$$
\aligned &\Lambda_{w} ={\Eu}(\een_{\Vb,w}^*\oplus\nen_{\Vb,w}),
\hfill\cr &\Lambda_{ws} ={\Eu}(\een_{\Vb,ws}^*\oplus\nen_{\Vb,ws}),
\hfill\cr &\Lambda^{s}_{w,w} =\Lambda^{s}_{w,ws}
=\Eu(\een_{\Vb,ws,w}^*\oplus\nen_{\Vb,w})^{-1}, \hfill\cr
&\Lambda^{s}_{ws,w} =\Lambda^{s}_{ws,w}
=\Eu(\een_{\Vb,w,ws}^*\oplus\nen_{\Vb,ws})^{-1}.
\endaligned
\leqno(2.13)
$$
We have also
$$\gathered
\Eu(\een_{\Vb,w})\,\Eu(\den_{\Vb,ws,w})=
\Eu(\een_{\Vb,ws})\,\Eu(\den_{\Vb,w,ws}),\hfill\cr
\Eu(\een_{\Vb,w,ws})\,\Eu(\den_{\Vb,w,ws})= \Eu(\een_{\Vb,w}).
\endgathered$$
Thus (2.13) yields
$$\gathered
\Lambda_w\,\Eu(\den_{\Vb,ws,w}^*)=
\Lambda_{ws}\,\Eu(\den_{\Vb,w,ws}^*), \hfill\cr \Lambda^{s}_{w,w}
=\Lambda^{s}_{w,ws} =\Lambda_{ws}^{-1}\,\Eu(\den_{\Vb,ws,w}^*)
=\Lambda_{w}^{-1}\,\Eu(\den_{\Vb,w,ws}^*) =\Lambda^{s}_{ws,ws}
=\Lambda^{s}_{ws,w}.
\endgathered
\leqno(2.14)$$ Therefore we get the following alternatives

\vskip1mm

\item{$\bullet$}
$\ib'=s(\ib)$, $w'=ws$ and $g_{w'}=
w's(f)\,\Eu(\den_{\Vb,ws,w}^*)\,\Lambda_{w'}^{-1}$. Thus (2.11)
holds with
$$g=s(f)(\chi_{l+1}-\chi_{l})^{h_{i_l,i_{l+1}}}.$$

\vskip1mm

\item{$\bullet$}
$\ib'=\ib$, $w'=w$ and $g_{w'}=
w'(f)\,\Eu(\den_{\Vb,w,ws}^*)\,\Lambda_{w'}^{-1}$. Thus (2.11) holds
with
$$g=f(\chi_{l+1}-\chi_{l})^{h_{i_{l+1},i_l}}.$$

%
%
%
%

\qed

\vskip3mm

\subhead 2.24.~Description of the grading on $\Zb_\Vb$ and ${}^\delta\Zb_\Vb$\endsubhead 
For each sequence $\ib=(i_1,i_2,\dots, i_m)$ and each integers $k$,
$l$ as above we abbreviate
$$\gathered
\varkappa_\ib(k)=\varkappa_{\ib,\ib}(k),\quad
\sigma_{\ib}(l)=\sigma_{s(\ib),\ib}(l),\cr
h_\ib(l)=h_{i_l,i_{l+1}}\
\roman{if}\ s_l(\ib)\neq\ib,\quad h_\ib(l)=-1\ \roman{if}\
s_l(\ib)=\ib,\cr
a_\ib(l)=h_\ib(l)+h_{s_l(\ib)}(l),\cr
a_\ib(l)=-i_l\cdot i_{l+1}\
\roman{if}\ s_l(\ib)\neq\ib,\quad 
a_\ib(l)=-2\
\roman{if}\ s_l(\ib)=\ib.
\endgathered\leqno(2.15)$$
The action of the elements $1_\ib$, $\varkappa_{\ib}(k)$, $\sigma_{\ib}(l)$ of $\Zb_\Vb$ on
$\Fb_\Vb$ yields linear operators in $\End(\Fb_\Vb)$. Let us denote
them by $\tilde 1_\ib$, $\tilde \varkappa_{\ib}(k)$,
$\tilde\sigma_{\ib}(l)$ respectively. Recall that $\Fb_\Vb$ is a
faithful left graded $\Zb_\Vb$-module by Section 1.6 and Lemma 1.8, and that the isomorphism in (2.4)
$$\Fb_\Vb=\bigoplus_\ib\kb[x_\ib(1),x_\ib(2),\dots x_\ib(m)]$$
is a graded $\kb$-vector space isomorphism.
Hence Proposition 2.23 and (2.4) imply the following.

\proclaim{2.25. Corollary} The grading on the $\k$-algebra $\Zb_\Vb$
is uniquely determined by the following rules. For each $\ib\in
I^\nu$ and each $k,l=1,2,\dots m$, $l\neq m$, we have :
$1_\ib$ has the degree 0,
$\varkappa_{\ib}(k)$ has the degree $2$,
and
$\sigma_{\ib}(l)$ has the degree $2h_\ib(l)$.
\endproclaim

\noindent By (1.11) we have
$${}^\delta\Zb_{s_l(\ib),\ib}=
\Zb_{s_l(\ib),\ib}[d_{s_l(\ib)}-d_{\ib}].$$ Further \cite{L1,
lem.~1.6$(c)$} yields
$$d_\ib=\ell_\nu+\sum_{l'\leqslant l}h_{i_{l},i_{l'}}$$
(note that our notation for flags is opposite to the one in
loc.~cit.). Thus we have
$$d_{s_l(\ib)}=d_\ib-h_{s_l(\ib)}(l)+h_\ib(l).$$
In particular the grading in ${}^\delta\Zb_\Vb$ obeys the following rules :
$1_\ib$ has the degree 0,
$\varkappa_{\ib}(k)$ has the degree 2,
and
$\sigma_{\ib}(l)$ has the degree $a_\ib(l)$.

\subhead 2.26.~The PBW theorem for $\Zb_\Vb$\endsubhead View $\Fb_\Vb$ as a graded (commutative)
$\k$-algebra as in Section 2.14. The graded $\Sb_\Vb$-algebra
$\Zb_\Vb$ is also a left graded $\Fb_\Vb$-module : we let
$f(x_\ib(1),\dots,x_\ib(m))\in\Fb_\ib$ act on $\Zb_\Vb$ as the
operator $$z\mapsto
f(\varkappa_{\ib}(1),\dots,\varkappa_{\ib}(m))\star z.$$
By Lemma 2.9$(b)$ the left graded $\Sb_\Vb$-submodule
$\Zb^{\leqslant x}_\Vb\subset\Zb_\Vb$ is indeed a left graded
$\Fb_\Vb$-submodule for each $x\in\Sen$.

\proclaim{2.27.~Lemma} We have $\Zb_\Vb^{\leqslant x}=
\bigoplus_{w\leqslant x}\Fb_\Vb\star[Z_\Vb^{w}]$ for each
$w\in\Sen$. In particular $\Zb_\Vb$ is a free left graded
$\Fb_\Vb$-module of rank $m!$.
\endproclaim

\noindent{\sl Proof :} The proof is the same as in \cite{CG,
sec.~7.6.11} which considers the equivariant $K$-theory of the
Steinberg variety. It uses a decomposition of $Z_\Vb^{\leqslant x}$
as a disjoint union of affine cells. Details are left to the reader.

\qed



\subhead 2.28.~Examples\endsubhead $(a)$ Set $I=\{i\}$
and $\nu=m\,i$. So $I^\nu=\{\ib\}$ with
$\ib=(i,i,\dots, i)$, $G_\Vb=\GL(\CC^m)$, $E_\Vb=\{0\}$, $\widetilde
F_\Vb=F_\Vb=F$ and $Z_\Vb=F\times F$. We have $\Fb_\Vb\simeq\Pb_\Vb$
as a left graded $\Sb_\Vb$-module. We have also
$\Zb_\Vb\simeq\End_{\Sb_\Vb}(\Fb_\Vb)$ as a graded
$\Sb_\Vb$-algebra. It is the nil-Hecke ring. Note that $\Fb_\Vb$ is
a projective left graded $\Zb_\Vb$-module such that
$\Zb_\Vb=\bigoplus_{w\in\Sen_m}\Fb_\Vb[2\ell(w)]$ as a left graded
$\Zb_\Vb$-module.

\vskip1mm

$(b)$ Set $I=\{i_1,i_2,\dots,i_m\}$, $H=\emptyset$ and
$\nu=i_1+i_2+\cdots +i_m$. So $I^\nu\simeq\Sen_m$ and
$F_\ib\simeq\{\bullet\}$ for all $\ib$. We have also
$G_\Vb\simeq(\CC^\times)^d$, $E_\Vb=\{0\}$, $\widetilde F_\Vb=F_\Vb$
and $Z_\Vb=F_\Vb\times F_\Vb.$ Thus $\Fb_\Vb$ is a free left graded
$\Sb_\Vb$-module of rank $m!$ and $\Zb_\Vb=\End_{\Sb_\Vb}(\Fb_\Vb)$
as a graded $\Sb_\Vb$-algebra.

\vskip1mm

$(c)$ Set $I=\{i,j\}$ with $j\to i$ and $\nu=i+j$. So
$I^\nu=\{\ib,\ib'\}$ with $\ib=(i,j)$, $\ib'=(j,i)$ and $F_\ib\simeq
F_{\ib'}\simeq\{\bullet\},$ $E_\Vb\simeq\CC,$ $\widetilde
F_\ib\simeq\CC$, $\widetilde F_{\ib'}\simeq\{\bullet\}$,
$Z_{\ib,\ib}\simeq\CC$ and $Z_{\ib,\ib'}\simeq Z_{\ib',\ib}\simeq
Z_{\ib',\ib'}\simeq \{\bullet\}.$ In particular $\Fb_\Vb$ is a free
left graded $\Sb_\Vb$-module of rank 2.

\vskip3mm

\subhead 3.~KLR-algebras\endsubhead

In this section we compute the Ext-algebras introduced in the first
section.

\subhead 3.1.~The algebra $\Rb_\Vb$ and its canonical representation on $\Fb_\Vb$\endsubhead 
Let $\nu$ be a non zero element of $\NN_mI$. Fix $\Vb\in\Vcb_\nu$. Let
$\Rb_\Vb$ be the KLR-algebra associated with $\Vb$. It is a graded
$\k$-algebra with $1$. By \cite{KL1, sec.~2.3} there is a faithful left
representation of $\Rb_\Vb$ on
$$\Qb_\Vb=\bigoplus_{\ib\in I^\nu}\Qb_\ib,\quad
\Qb_\ib=\k[x_\ib(1),\dots x_\ib(m)].$$ Recall that (2.4) yields a
$\kb$-vector space isomorphism
$$\Qb_\Vb=\Fb_\Vb.\leqno(3.1)$$
From now on we'll write $\Fb_\Vb$ everywhere instead of $\Qb_\Vb$.
Next, Proposition 2.23 and Section 2.24 yield explicit linear
operators
$$\tilde 1_\ib,\,\tilde \varkappa_{\ib}(k),\,
\tilde\sigma_{\ib}(l)\in\End (\Fb_\Vb),\quad\ib\in I^\nu,\,
1\leqslant k,l\leqslant m,\,l\neq m.$$ They are given as follows

\vskip1mm \item{$\bullet$} $\tilde 1_\ib$ is the projection on
$\Fb_\ib$ relatively to $\bigoplus_{\ib'\neq\ib}\Fb_{\ib'}$,

\vskip1mm \item{$\bullet$} $\tilde\varkappa_{\ib}(k)$ acts by zero
on $\Fb_{\ib'}$, $\ib'\neq\ib$, and by multiplication by $x_\ib(k)$
on $\Fb_\ib$,

\vskip1mm \item{$\bullet$} $\tilde\sigma_{\ib}(l)$ acts by zero on
$\Fb_{\ib'}$, $\ib'\neq\ib$, and by the following rules \vskip-2mm

$$\aligned
&f\mapsto (x_\ib(l+1)-x_\ib(l))^{-1}(s_l(f)-f) \ \roman{if}\
s_l(\ib)=\ib, \hfill\cr &f\mapsto
(x_{s_l(\ib)}(l+1)-x_{s_l(\ib)}(l))^{h_\ib(l)}s_l(f)\
\roman{if}\ s_l(\ib)\neq \ib.
\endaligned\leqno(3.2)$$

\noindent By loc.~cit.~ the embedding $\Rb_\Vb\subset\End(\Fb_\Vb)$
identifies $\Rb_\Vb$ with the $\k$-subalgebra generated by $\tilde
1_\ib$, $\tilde\varkappa_{\ib}(k)$ and $\tilde\sigma_{\ib}(l)$ for
all $\ib$, $k$, $l$ as above.

\subhead 3.2.~The grading on $\Rb_\Vb$ and $\Fb_\Vb$\endsubhead The grading on $\Rb_\Vb$ is given by the
following rules, see \cite{KL1} :
$\tilde 1_\ib$ has the degree 0,
$\tilde \varkappa_{\ib}(k)$ has the degree 2,
and
$\tilde \sigma_{\ib}(l)$ has the degree $a_\ib(l)$.
Given a sequence $\ib_0\in I^\nu$ there is a unique grading on
$\Fb_\Vb$ such that $\Fb_\Vb$ is a graded left $\Rb_\Vb$-module and
the unit of $\Fb_{\ib_0}$ has the degree 0. Compare \cite{KL1,
sec.~2.3}. This grading depends on the choice of the sequence
$\ib_0$. By Corollary 2.25 there is a canonical graded $\kb$-vector
space isomorphism
$$\Fb_\Vb={}^\delta\Fb_\Vb,\leqno(3.3)$$ up to a shift by an integer
that we'll ignore from now on. Here the lhs denotes the underlying graded $\kb$-vector space of the graded
$\Rb_\Vb$-module $\Qb_\Vb$ while the rhs is the underlying $\kb$-vector space of the graded ${}^\delta\Zb_\Vb$-module
from Section 1.9.

\subhead 3.3.~Remarks\endsubhead $(a)$ The formulas in Proposition
2.23 differ indeed from the ones in \cite{KL1, sec.~2.3} by a sign.
However it is explained in \cite{KL2} how to modify the algebra
$\Rb_\Vb$ so that all the results of \cite{KL1} remain true and so
that this difference of sign disappear. From now on we'll ignore
this point, and we'll refer only to \cite{KL1} to simplify.

$(b)$ In \cite{KL1} the authors assume that two different vertices
in $\Gamma$ may be joined by at most one arrow. The quivers and the
algebras we consider are more general. They appear in \cite{R,
sec.~3.2.4}. Since this generalization does not affect the rest of
the paper, from now on we'll assume that \cite{KL1} works in this
greater generality. Indeed, we'll only use the fact that  $\Rb_\Vb$
admits a faithful representation on $\Fb_\Vb$ given by the same
formulas as in Section 3.1 and the PBW-theorem. Both facts are
proved for general quivers in \cite{R, prop.~3.12, thm.~3.7}. Note
that some algebra are also defined in \cite{KL2} for arbitrary
quivers, but they are different from the ones considered here and in
\cite{R}.

$(c)$ The $\kb$-algebra $\Rb_\Vb$ has the following presentation,
see \cite{R, def.~3.2.1} (we'll not use this) : it is generated by
$1_\ib$, $\varkappa_\ib(k)$, $\tau_\ib(l)$ with $\ib\in I^\nu$,
$1\leqslant k,l\leqslant m$ and $l\neq m$, modulo the following
defining relations

\vskip1mm

\item{$\bullet$}
$1_\ib\,1_{\ib'}=\delta_{\ib,\ib'}1_\ib$,

\vskip1mm

\item{$\bullet$}
$\tau_\ib(l)=1_{s_l(\ib)}\tau_\ib(l)1_\ib$,

\vskip1mm

\item{$\bullet$}
$\varkappa_\ib(k)=1_{\ib}\varkappa_\ib(k)1_\ib$,

\vskip1mm

\item{$\bullet$}
$\varkappa_\ib(k)\varkappa_{\ib'}(k')=\varkappa_{\ib'}(k')\varkappa_\ib(k)$,

\vskip1mm

\item{$\bullet$}
$\tau_{s_l(\ib)}(l)\tau_{\ib}(l)=
Q_{\ib,l}(\varkappa_\ib(l),\varkappa_\ib(l+1))$,

\vskip1mm

\item{$\bullet$}
$\tau_{s_l(\ib)}(l')\tau_{\ib}(l)=\tau_{s_{l'}(\ib)}(l)\tau_\ib(l')$
if $|l-l'|>1$,

\vskip1mm

\item{$\bullet$}
$\tau_{s_ls_{l+1}(\ib)}(l+1)\tau_{s_{l+1}(\ib)}(l)\tau_{\ib}(l+1)-
\tau_{s_{l+1}s_l(\ib)}(l)\tau_{s_{l}(\ib)}(l+1)\tau_{\ib}(l)=
(Q_{\ib,l}(\varkappa_\ib(l+2),\varkappa_\ib(l+1))
-Q_{\ib,l}(\varkappa_\ib(l),\varkappa_\ib(l+1)))
(\varkappa_\ib(l+2)-\varkappa_\ib(l))^{-1}$ if $s_{l,l+2}(\ib)=\ib$
and 0 else,

\vskip1mm

\item{$\bullet$}
$\tau_{\ib}(l)\varkappa_{\ib}(k)-\varkappa_{s_{l}(\ib)}(s_l(k))\tau_\ib(l)
=\cases
-1_\ib&\ \roman{if}\  k=l,\, s_l(\ib)=\ib, \cr
1_\ib&\ \roman{if}\  k=l+1,\, s_l(\ib)=\ib,\cr
0&\ \roman{else}.
\endcases$

\vskip1mm

\noindent Here we have set $ s_{l,l+2}=s_ls_{l+1}s_l$ if $l\neq
m-1,m$ and
$$Q_{\ib,l}(u,v)=\cases(-1)^{h_\ib(l)}(u-v)^{a_\ib(l)}&\ \roman{if}\
s_l(\ib)\neq \ib, \cr 0&\ \roman{else}.\endcases$$ The element
$\tau_\ib(l)$ acts on $\Fb_\Vb$ as the operator
$(-1)^{h_\ib(l)}\tilde\sigma_\ib(l),$ while the element
$\varkappa_\ib(k)$ acts by multiplication by $x_\ib(k)$.

\subhead 3.4.~The PBW Theorem for $\Rb_\Vb$ and the isomorphism
$\Rb_\Vb={}^\delta\Zb_\Vb$\endsubhead Note that $\Rb_\Vb$ is a free left graded
$\Fb_\Vb$-module such that $f(x_\ib(1),\dots,x_\ib(m))$ acts on
$\Rb_\Vb$ by the left multiplication with the element
$f(\tilde\varkappa_{\ib}(1),\dots,\tilde\varkappa_{\ib}(m))$. Let
the symbol $\star$ denote the left $\Fb_\Vb$-action on $\Rb_\Vb$.
One constructs a $\Fb_\Vb$-basis of $\Rb_\Vb$ in the following way.
For each permutation $w\in\Sen_m$ we choose a reduced decomposition
$w=s_{l_1}s_{l_2}\cdots s_{l_r}$ with $r$ an integer $\geqslant 0$
and $l_1,l_2,\dots, l_r\in\{1,2,\dots m-1\}$. For each $\ib\in
I^\nu$ let $\tilde\sigma_\ib(w)\in \tilde 1_\ib\star\Rb_\Vb$ be
given by $\tilde\sigma_\ib(w)=1_\ib$ if $r=0$, and
$$\tilde\sigma_\ib(w)= \tilde\sigma_{s_{l_1}(\ib)}(l_1) \star
\tilde\sigma_{s_{l_2}s_{l_1}(\ib)}(l_2) \star\cdots\star
\tilde\sigma_{w^{-1}(\ib)}(l_r)\ \roman{else}.$$ Observe that
$\tilde\sigma(w)=\sum_\ib\tilde\sigma_\ib(w)$ may depend on the
choice of the reduced decomposition of $w$.
The following is proved in \cite{KL1, thm.~2.5}.

\proclaim{3.5. Proposition} We have
$\Rb_\Vb=\bigoplus_{w\in\Sen_m}\Fb_\Vb\star\tilde\sigma(w)$ as a
left $\Fb_\Vb$-module.
\endproclaim

Restricting the $\Fb_\Vb$-action on $\Rb_\Vb$ to the subalgebra
$\Sb_\Vb\subset\Fb_\Vb$ we get a structure of graded
$\Sb_\Vb$-algebra on $\Rb_\Vb$. See Sections 2.13-2.15 for details.
The first main result of this paper is the following.

\proclaim{3.6. Theorem} There is an unique graded $\Sb_\Vb$-algebra
isomorphism
$$\Psi_\Vb:\Rb_\Vb\to{}^\delta\Zb_\Vb$$ such that (3.3) intertwines the
actions of $\Rb_\Vb$, ${}^\delta\Zb_\Vb$.
\endproclaim

\noindent{\sl Proof :} First, recall that ${}^\delta\Fb_\Vb$ is a
faithful left graded ${}^\delta\Zb_\Vb$-module and that $\Rb_\Vb$ is
the graded $\k$-subalgebra of $\End(\Fb_\Vb)$ generated by the
operators $\tilde 1_\ib$, $\tilde\varkappa_{\ib}(k)$ and
$\tilde\sigma_{\ib}(l)$. Thus there is a unique injective graded
$\k$-algebra homomorphism
$$\Psi_\Vb:\Rb_\Vb\to{}^\delta\Zb_\Vb,\quad
\Psi_\Vb(\tilde 1_\ib)= 1_\ib, \quad
\Psi_\Vb(\tilde\varkappa_{\ib}(k))= \varkappa_{\ib}(k), \quad
\Psi_\Vb(\tilde\sigma_{\ib}(l))=\sigma_{\ib}(l).$$ We must prove
that $\Psi_\Vb$ is a surjective map. The map $\Psi_\Vb$ is a left
graded $\Fb_\Vb$-module homomorphism by Proposition 2.23, Corollary
2.25 and Sections 3.1, 3.2. For each $x\in\Sen$ we set
$$\Rb_\Vb^{\leqslant x}=\bigoplus_{w\leqslant
x}\Fb_\Vb\star\tilde\sigma(w),$$ a left graded $\Fb_\Vb$-submodule
of $\Rb_\Vb$. We abbreviate $\Rb_\Vb^{e}=\Rb_\Vb^{\leqslant e}.$ The
proof of the theorem consists of two steps. First we prove that
$\Psi_\Vb(\Rb_\Vb^{\leqslant x})\subset\Zb_\Vb^{\leqslant x}$. Then
we prove that this inclusion is an equality.

\vskip1mm

{\sl Step 1 :} Since $\Psi_\Vb$ is a left $\Fb_\Vb$-module
homomorphism it is enough to prove that
$\Psi_\Vb(\tilde\sigma(w))\in\Zb_\Vb^{\leqslant w}$ for each $w$. By
Lemma 2.9$(b)$ and an easy induction on the length of $w$ it is
enough to prove the following
$$\Psi_\Vb(1)\subset\Zb_\Vb^{e},\quad
\Psi_\Vb(\tilde\sigma(s_l))\subset\Zb_\Vb^{\leqslant
s_l},\quad\forall l.$$ This is obvious, because $1\in\Zb_\Vb^e$ and
$\sigma_{\ib}(l)\in\Zb_\Vb^{\leqslant s_l}$ for each $\ib$, see
Section 2.22.

\vskip1mm

{\sl Step 2 :} Lemma 2.27 implies that $\Zb^e_\Vb$ is the free
$\Fb_\Vb$-module of rank one generated by $[Z^e_\Vb]$. Therefore we
have
$$\Psi_\Vb(\Rb_\Vb^{e})=
\Psi_\Vb(\Fb_\Vb\star 1)= \Fb_\Vb\star\Psi_\Vb(1)=
\Fb_\Vb\star[Z_\Vb^e]= \Zb_\Vb^{e}.$$ Next, fix an integer
$l=1,2,\dots m-1$. We claim that $$\Psi_\Vb(\Rb_\Vb^{\leqslant
s_l})=\Zb_\Vb^{\leqslant s_l}.$$ By Lemma 2.27, Proposition 3.5 we
have
$$\aligned
&\Rb_\Vb^{\leqslant s_l}= (\Fb_\Vb\star
1)\oplus(\Fb_\Vb\star\tilde\sigma(s_l)),\hfill\cr
&\Zb_\Vb^{\leqslant s_l}=
(\Fb_\Vb\star[Z^e_\Vb])\oplus(\Fb_\Vb\star[Z^{s_l}_\Vb]).
\endaligned$$ Thus it is enough to observe that
$$\Psi_\Vb(1)=[Z^{e}_\Vb],\quad\Psi_\Vb(\tilde\sigma(l))=[Z^{s_l}_\Vb].$$
See the definition of $\sigma(l)$ in Section 2.22 for the second
identity. To complete the proof of Step 2 we are reduced to prove
the following.

\proclaim{3.7. Lemma} If $\ell(s_lw)=\ell(w)+1$ we have
$[Z^{s_l}_\Vb]\star[Z^w_\Vb]=[Z^{s_lw}_\Vb]$ in $\Zb_\Vb^{\leqslant
s_lw}/\Zb_\Vb^{<s_lw}$.
\endproclaim

\noindent{\sl Proof :} Assume that $\ell(s_lw)=\ell(w)+1$. We'll
abbreviate $s=s_l$ as above. By Lemma 2.9$(b)$ and Lemma 2.27 there
is an unique element $c\in\Fb_\Vb$ such that
$$[Z^{s}_\Vb]\star[Z^w_\Vb]=c\star[Z^{sw}_\Vb]\ \roman{in}\ \Zb_\Vb^{\leqslant
s_lw}/\Zb_\Vb^{<s_lw}.$$ We must prove that $c=1$.

For each $x\in\Sen$ we abbreviate $[Z^x_\Vb]=[Z^x_\Vb]\otimes 1$, an
element of $\Hc_*(Z_\Vb,\k)$. For each $y,z\in\Sen$ there is a
unique element $\Lambda^x_{y,z}\in\Kb_\Vb$ such that
$$[Z^x_\Vb]=\sum_{y,z}\Lambda^x_{y,z}\psi_{y,z}.$$
Compare (2.9). Since $\phi_{\Vb,y,yx}$ is a smooth point of
$Z^x_\Vb$ we have also
$$\Lambda^x_{y,yx}=\Eu(Z^x_\Vb,\phi_{\Vb,y,yx})^{-1}.$$
Hence, in the expansion of the element $[Z^{sw}_\Vb]$ in the
$\Kb_\Vb$-basis $(\psi_{y,z})$  the coordinate along the vector
$\psi_{x,xsw}$ is equal to
$$\Lambda_{x,xsw}^{sw}=\Eu(Z^{sw}_\Vb,\phi_{\Vb,x,xsw})^{-1}.$$
On the other hand, since $\Lambda^w_{x,xsw}=0$ and
$$[Z^s_\Vb]=\sum_{x}(\Lambda^s_{x,x}\psi_{x,x}+\Lambda^s_{x,xs}\psi_{x,xs}),$$
the coordinate of $[Z^{s}_\Vb]\star[Z^w_\Vb]$ along $\psi_{x,xsw}$
is equal to
$$\Lambda_{x,xs}^{s}\Lambda^w_{xs,xsw}\Lambda_{xs}=
\Eu(Z^{s}_\Vb,\phi_{\Vb,x,xs})^{-1}
\Eu(Z^{w}_\Vb,\phi_{\Vb,xs,xsw})^{-1}\Lambda_{xs}.$$ Thus we must
check that
$$\Eu(Z^{s}_\Vb,\phi_{\Vb,x,xs})\,
\Eu(Z^{w}_\Vb,\phi_{\Vb,xs,xsw})= \Eu(Z^{sw}_\Vb,\phi_{\Vb,x,xsw})\,
\Lambda_{xs}.$$ This follows from Lemma 3.8 below.

\qed

\proclaim{3.8. Lemma} (a) For each $x,y\in\Sen$ we have
$$\aligned
\eu(O^y_\Vb,\phi_{\Vb,x,xy})&=\eu(\nen_{\Vb,x}\oplus\men_{\Vb,xy,x}),
\hfill\cr \Eu(Z^{y}_\Vb,\phi_{\Vb,x,xy})&=
\Eu(O^{y}_\Vb,\phi_{\Vb,x,xy})\,\Eu(\een_{\Vb,x,xy}^*), \hfill\cr
\Lambda_x=\Eu(Z^{e}_\Vb,\phi_{\Vb,x,x})
&=\Eu(F_\Vb,\phi_{\Vb,x})\,\Eu(\een_{\Vb,x}^*).
\endaligned$$

(b) For each $w,x,y\in\Sen$ such that $\ell(xy)=\ell(x)+\ell(y)$ we
have
$$\aligned
\eu(O_\Vb^{xy},\phi_{\Vb,w,wxy})\,\eu(F_\Vb,\phi_{\Vb,wx})&=
\eu(O^x_\Vb,\phi_{\Vb,w,wx})\,\eu(O_\Vb^{y},\phi_{\Vb,wx,wxy}),
\hfill\cr \Eu(\een^*_{\Vb,w,wxy}\oplus \een^*_{\Vb,wx})&=
\Eu(\een^*_{\Vb,w,wx}\oplus \een^*_{\Vb,wx,wxy}).
\endaligned$$
\endproclaim

\noindent{\sl Proof :} Part $(a)$ is left to the reader. Let us
prove $(b)$. Set $$\Delta(y)^-=y(\Delta^+)\cap\Delta^-,
\quad\Delta(y)^+=y(\Delta^-)\cap\Delta^+.$$ For each $x,y\in\Sen$
the $T_\Vb$-module $\men_{\Vb,xy,x}$ is the sum of the root
subspaces whose weight belong to the set
$x(\Delta(y)^-)\cap\Delta_\Vb$. Thus, by $(a)$, the first claim is
equivalent to the following equality
$$w(\Delta(xy)^-)\cap\Delta_\Vb=
w\bigl(\Delta(x)^-\sqcup x(\Delta(y)^-)\bigr)\cap\Delta_\Vb.
$$
This equality is a consequence of the following well-known formula
$$\ell(xy)=\ell(x)+\ell(y)\Rightarrow
\Delta(xy)^-=\Delta(x)^-\sqcup x(\Delta(y)^-).$$

Now, let $\Xi_\Vb\subset\Delta$ be the set of weights of the
$T_\Vb$-module $E_\Vb$. Note that a weight subspace appear in
$E_\Vb$ with the multiplicity at most one. So the character of the
$T_\Vb$-module $\een_{\Vb,x,xy}$ is the sum of the roots $\a$ which
belong to the set
$$x(\Delta^+\setminus\Delta(y)^+)\cap\Xi_\Vb.$$
Let $$A=(\Delta^+\setminus\Delta(xy)^+)\sqcup x(\Delta^+),\quad
B=(\Delta^+\setminus\Delta(x)^+)\sqcup
x(\Delta^+\setminus\Delta(y)^+).$$ The character of the
$T_\Vb$-module $\een^*_{\Vb,w,wxy}\oplus \een^*_{\Vb,wx}$ is
$\sum\a$ where $\a$ runs over the set $w(A)\cap\Xi_\Vb$. The
character of the $T_\Vb$-module $\een^*_{\Vb,w,wx}\oplus
\een^*_{\Vb,wx,wxy}$ is $\sum\b$ where $\b$ runs over the set
$w(B)\cap\Xi_\Vb$. Since $\ell(xy)=\ell(x)+\ell(y)$ we have
$$\Delta(xy)^+=\Delta(x)^+\sqcup x(\Delta(y)^+)$$
This implies that $A=B$. The second claim follows.

\qed

\vskip3mm

\subhead 3.9\endsubhead From now on we'll use freely the isomorphism
$\Psi_\Vb$ to identify the graded $\Sb_\Vb$-algebras $\Rb_\Vb$ and
${}^\delta\Zb_\Vb$. In particular we'll abbreviate $\tilde 1_\ib=1_\ib$,
$\tilde\varkappa_{\ib}(k)=\varkappa_\ib(k)$ and $\tilde
\sigma_{\ib}(l)=\sigma_{\ib}(l)$. When no confusion is possible
we'll also identify the graded $\k$-vector spaces $\Fb_\Vb$ and
${}^\delta\Fb_\Vb$ as in (3.3).

\subhead 4. Canonical bases and projective graded $\Rb_\Vb$-modules
\endsubhead

In this section, using the computations of Section 3 we prove a
conjecture of \cite{KL1} yielding a categorification of the
canonical basis of the negative part of Drinfeld-Jimbo quantized
enveloping algebra.

\subhead 4.1.~The Grothendieck group of $\Rb_\Vb$\endsubhead This subsection is reminder on the
Grothendieck group of the graded ring $\Rb_\Vb$. See \cite{KL1} for
details. It does not contain new results.
The graded $\kb$-algebra $\Rb_\Vb$ is finite dimensional over its center
$\Sb_\Vb$, a commutative graded $\kb$-subalgebra, and the quotient
$\Rb_\Vb/\Sb_\Vb^+\Rb_\Vb$ is finite-dimensional. Therefore any simple
object of $\Rb_\Vb$-$\modb$ is finite-dimensional.
Further there is a finite number of simple modules in
$\Rb_\Vb$-$\modb^f$.
Thus $G(\Rb_\Vb)$ is a free Abelian group of finite rank with
a basis formed by the classes of the
simple objects of $\Rb_\Vb$-$\modb$, while
$K(\Rb_\Vb)$ is a free Abelian group with
a basis formed by the classes of the
indecomposable projective objects. Both are free $\Ac$-modules where
$q$ shifts the grading by $1$. For each $\Vb,\Vb'\in\Vcb_\nu$ there is a canonical isomorphism
$\Rb_\Vb=\Rb_{\Vb'}$. Thus, from now on we'll abbreviate
$\Rb_\nu=\Rb_\Vb$.
Consider the Abelian group
$$K(\Rb)=\bigoplus_\nu K(\Rb_\nu).$$ 
Recall that $K(\Rb)$ is a $\Ac$-module by (0.2).
The $\Ac$-action on $K(\Rb)$ is the same as the one in \cite{KL1, sec.~2.5},
because both (0.2) and the grading on $\Rb_\nu$ are opposite to the
conventions in loc.~cit. We equip $K(\Rb)$ with an associative unital
$\Ac$-algebra structure as follows. Fix $\nu,\nu'\in\NN I$.
Set $\nu''=\nu+\nu'$ and
$m=|\nu|$. Given sequences $\ib\in I^\nu$,
$\ib'\in I^{\nu'}$ we write $\ib''=\ib\ib'$, see the notation in (0.1).
There is an unique
inclusion of graded $\k$-algebras
$$\Rb_\nu\otimes\Rb_{\nu'}\to\Rb_{\nu''}\leqno(4.1)$$ such that
for each $\ib$, $\ib'$, $k$, $l$ we have
$$1_\ib\otimes 1_{\ib'}\mapsto1_{\ib''}, \quad
\varkappa_{\ib}(k)\otimes 1_{\ib'}\mapsto \varkappa_{\ib''}(k),\quad
1_\ib\otimes \varkappa_{\ib'}(k)\mapsto \varkappa_{\ib''}(m+k),$$
$$\sigma_{\ib}(l)\otimes 1_{\ib'}\mapsto \sigma_{\ib''}(l),
\quad 1_\ib\otimes \sigma_{\ib'}(l)\mapsto \sigma_{\ib''}(m+l).$$
Let $1_{\nu,\nu'}$ be the image of the identity element by the map (4.1).
The induction yields an additive functor
$$\Ind_{\nu,\nu'}:
\left\{\gathered
\Rb_\nu\text{-}\modb\times\Rb_{\nu'}\text{-}\modb\to
\Rb_{\nu''}\text{-}\modb,\cr (M,M')\mapsto
\Rb_{\nu''}1_{\nu,\nu'}\otimes_{\Rb_\nu\otimes\Rb_{\nu'}}(M\otimes
M').
\endgathered\right.$$
It takes projectives to projectives and
it commutes with the shift of the grading. Thus it yields an
$\Ac$-linear homomorphism
$$
K(\Rb_\nu)\otimes_\Ac K(\Rb_{\nu'})\to K(\Rb_{\nu''}).$$ Taking the
sum over all $\nu,\nu'$ we get an associative unital
$\Ac$-algebra structure on $K(\Rb)$.

\subhead 4.2.~The projective $\Rb_\nu$-module $\Rb_\yb$\endsubhead Assume that $\nu\in\NN I$.
Given a pair $\ibt=(\ib,\ab)\in Y_\nu$ we define
an object $\Rb_{\ibt}$ in $\Rb_\nu$-$\projb$ as follows, see
\cite{KL1, sec.~2.5}.

\vskip1mm

\item{$\bullet$}
If $I=\{i\}$, $\nu=m\,i$, $\ib=i^m$ and $\ab=m$ then $\yb=(i,m)$ and
we set
$$\Rb_{\ibt}=\Rb_{i,m}=\Fb_\Vb[\ell_m].$$ As a left graded $\Rb_\nu$-modules we
have
$\Rb_{\nu}\simeq\bigoplus_{w\in\Sen_m}\Rb_{i,m}[2\ell(w)-\ell_m]$,
see Example 2.28$(a)$. So $\Rb_{i,m}$ is a direct summand of
$\Rb_\nu[\ell_m]$. We choose once for all an idempotent
$1_{i,m}\in\Rb_\nu$ such that
$$\Rb_{i,m}=(\Rb_\nu\star 1_{i,m})[\ell_m].$$

\item{$\bullet$}
If $\ib=(i_1,\dots, i_k)$ and $\ab=(a_1,\dots a_k)$ we define the
idempotent $1_{\ibt}\in\Rb_\nu$ as the image of the element
$\bigotimes_{l=1}^k1_{i_l,a_l}$ by the inclusion of graded
$\k$-algebras
$\bigotimes_{l=1}^k\Rb_{a_li_l}\!\!\subset\Rb_{\nu}$ in (4.1).
Then we set $$\Rb_{\ibt}=(\Rb_\nu\star 1_{\ibt})[\ell_\ab].$$

\subhead 4.3\endsubhead
The $\Rb_\nu$-module $\Rb_{\ibt}$ satisfies the following
properties, see loc.~cit.~for details.

\vskip1mm

\item{$\bullet$}
Let $\ib'\in I^\nu$ be the sequence obtained by expanding the pair
$\ibt=(\ib,\ab)$. We have the following formula in
$\Rb_\nu$-$\projb$
$$\Rb_{\ib'}\simeq\bigoplus_{w\in\Sen_\ab}\Rb_{\ibt}[2\ell(w)-\ell_\ab]
=:[\ab]!\Rb_\yb.$$
Thus we have the following formula in $K(\Rb)$
$$[\Rb_{\ib'}]=[\ab]![\Rb_{\ibt}].\leqno{(4.2)}$$

\item{$\bullet$}
If $I=\{i\}$, $\nu=m\,i$, $\ib=i^m$ then
$$\Rb_{i^m}=\Rb_{m\,i}=[m]!\Rb_{i,m},\quad 1_{i^m}=1.$$

\item{$\bullet$}
Given $\ibt=(\ib,\ab)\in  Y_\nu$, $\ibt'=(\ib',\ab')\in Y_{\nu'}$ we
set $\ibt\ibt'=(\ib\ib',\ab\ab')$. We have an isomorphism of left
graded $\Rb_{\nu''}$-modules
$$\Ind_{\nu,\nu'}(\Rb_{\ibt}\otimes \Rb_{\ibt'})\simeq\Rb_{\ibt''}.
\leqno(4.3)$$

\vskip3mm

\subhead 4.4.~The quantum group $\fb$\endsubhead Set $\Kc=\QQ(q)$. Let $\fb$ be the negative
half of the quantum universal enveloping algebra associated with the
quiver $(I,H)$. It is the $\Kc$-algebra generated by elements
$\theta_i$, $i\in I$, with the defining relations
$$\sum_{a+b=1-i\cdot j}(-1)^a\theta_i^{(a)}\theta_j\theta_i^{(b)}=0,
\quad \forall i\neq j.$$
Here we have set $\theta_i^{(a)}=\theta_i^a/[a]!$ for each integer
$a> 0$. Let ${}_\Ac\fb\subset\fb$ be the $\Ac$-lattice generated by
all products of the elements $\theta_i^{(a)}$. We have the weight
decomposition
$${}_\Ac\fb=\bigoplus_{\nu\in\NN I}{}_\Ac\fb_{\nu},\quad
\fb=\bigoplus_{\nu\in\NN I}\fb_{\nu}.$$ For each pair
$\ibt=(\ib,\ab)\in Y_\nu$ with $\ib=(i_1,\dots, i_k)$,
$\ab=(a_1,\dots a_k)$ we write
$$\theta_{\ibt}=
\theta_{i_1}^{(a_1)}\theta_{i_2}^{(a_2)}\cdots\theta_{i_k}^{(a_k)}\in
{}_\Ac\fb_\nu.$$ Let $\Bc$ be the canonical basis of ${}_\Ac\fb$. We
write
$$\Bc_\ZZ=\{q^d\bb;\bb\in\Bc,d\in\ZZ\},\quad
\Bc_\nu=\Bc\cap{}_\Ac\fb_\nu.$$ By \cite{KL1, prop.~3.4, sec.~3.2}
there is an unique $\Ac$-algebra isomorphism
$$\g_\Ac:{}_\Ac\fb\to K(\Rb),\quad
\theta_{\ibt}\mapsto[\Rb_{\ibt}],\quad\ibt\in Y_\nu.\leqno(4.4)$$
The following is the second main result of the paper.

\proclaim{4.5.~Theorem} The map $\g_\Ac$ takes $\Bc_\ZZ$ to the
$\ZZ$-basis of $K(\Rb)$ consisting of the indecomposable projective
objects.
\endproclaim

\subhead 4.6.~Definition of the canonical basis $\Bc$\endsubhead Before the proof let us recall the
construction of $\Bc$ in \cite{L2}. Fix a non zero element $\nu\in\NN I$ and fix
$\Vb\in\Vcb_\nu$. Let $\Pcb_\Vb$ the set of isomorphism classes of
simple perverse sheaves $\Lc$ on $E_\Vb$ such that $\Lc[r]$ appears
as a direct summand of $\Lc_{\ib}$ for some $\ib\in I^\nu$,
$r\in\ZZ$. Such a perverse sheaf belongs to $\Dcb_{G_\Vb}(E_\Vb)$.
Let $\Qcb_\Vb$ the full subcategory of $\Dcb_{G_\Vb}(E_\Vb)$
consisting of all complexes that are isomorphic to finite direct
sums of complexes of the form $\Lc[r]$ for various $r\in\ZZ$ and
$\Lc\in\Pcb_\Vb$.

Let $K(\Qcb_\Vb)$ be the Abelian group with one generator $[\Lc]$
for each isomorphism class of objects of $\Qcb_\Vb$ and with
relations $[\Lc]+[\Lc']=[\Lc'']$ whenever $\Lc''$ is isomorphic to
$\Lc\oplus\Lc'$. It is a free $\Ac$-module such that $q\Lc=\Lc[1]$
and $q^{-1}\Lc=\Lc[-1]$, where $\Lc$ runs over $\Qcb_\Vb$. An
isomorphism $\Vb\simeq\Vb'$ in $\Vcb$ induces a canonical
isomorphism $K(\Qcb_\Vb)\simeq K(\Qcb_{\Vb'})$. Taking the direct
limit over the groupoid consisting of the objects of $\Vcb$ with their isomorphisms, we get
$$K(\Qcb)=\ind_\Vb K(\Qcb_\Vb).$$
Given $\nu,\nu'\in\NN I$ and $\Vb\in\Vcb_\nu$, $\Vb'\in\Vcb_{\nu'}$
we set $\nu''=\nu+\nu'$, $\Vb''=\Vb\oplus\Vb'$. Let
$$*:\Qcb_\Vb\times\Qcb_{\Vb'}\to \Qcb_{\Vb''}$$ be
Lusztig's induction functor \cite{L2, sec.~9.2.5}. There is an
unique associative $\Ac$-bilinear multiplication $\circledast$ on
$K(\Qcb)$ such that for each $\Lc\in\Qcb_\Vb$, $\Lc'\in\Qcb_{\Vb'}$
we have
$$\Lc\circledast\Lc'=(\Lc*\Lc')[m_{\nu,\nu'}],\quad
m_{\nu,\nu'}=\sum_{h\in H}\nu_{h'}\nu'_{h''}+\sum_{i\in
I}\nu_i\nu'_i.$$ There is also an unique $\Ac$-algebra isomorphism
$$\l_\Ac:K(\Qcb)\to{}_\Ac\fb,\quad
[{}^\delta\!\Lc_{\ibt}]\mapsto\theta_{\ibt},\quad \forall\ibt.
\leqno(4.5)$$ See \cite{L2, thm.~13.2.11}. The classes in
$K(\Qcb_\Vb)$ of the perverse sheaves of $\Pcb_\Vb$ form a
$\Ac$-basis of $K(\Qcb_\Vb)$. We have
$$\Bc_\nu=\{\bb_\Lc; \Lc\in\Pcb_\Vb\},\quad
\quad\bb_\Lc=\l_\Ac([\Lc]).$$ For any $\nu,\nu'\in\NN I$ and
$\ibt\in  Y_\nu$, $\ibt'\in  Y_{\nu'}$ we have  \cite{L2,
sec.~9.2.6-7}
$${}^\delta\!\Lc_{\ibt}\circledast{}^\delta\!\Lc_{\ibt'}
\simeq{}^\delta\!\Lc_{\ibt\ibt'}.\leqno(4.6)$$ If $\ib'\in I^\nu$ is
the expansion of a pair $\ibt=(\ib,\ab)\in Y_\nu$ then we have also
$$[{}^\delta\!\Lc_{\ib'}]=[\ab]![{}^\delta\!\Lc_{\ibt}].\leqno(4.7)$$

\subhead 4.7\endsubhead {\sl Proof of Theorem 4.5 :} Assume that
$\nu\in\NN I$ and $\Vb\in\Vcb_\nu$. Given $\Lc\in\Qcb_\Vb$ we
consider the left graded ${}^\delta\Zb_\Vb$-module given by
$$\Yb_\Lc=\Ext^*_{G_\Vb}({}^\delta\!\Lc_\Vb,\Lc).$$
We'll view it as a left graded $\Rb_\nu$-module via the isomorphism
$\Psi_\Vb$. For each complexes $\Lc,\Lc'\in\Qcb_\Vb$ and each
integer $d$ we have canonical isomorphisms
$$\Yb_{\Lc\oplus\Lc'}=\Yb_\Lc\oplus \Yb_{\Lc'},\quad
\Yb_{\Lc[r]}=\Yb_\Lc[r].$$ If $\Lc\in \Pcb_\Vb$ there is a sequence
$\ib\in I^\nu$ and an integer $r$ such that $\Lc[r]$ is a direct
summand of the semisimple complex ${}^\delta\!\Lc_\ib$. Further,  we
have
$$\Yb_{{}^\delta\!\Lc_{\ib}}=
\Ext^*_{G_\Vb}({}^\delta\!\Lc_\Vb,{}^\delta\!\Lc_{\ib})=
{}^\delta\Zb_\Vb\star 1_\ib\simeq\Rb_\ib\leqno(4.8)$$ as left graded
$\Rb_\nu$-modules. Therefore $\Yb_\Lc$ belongs to
$\Rb_\nu$-$\projb$. Since any element in $\Qcb_\Vb$ is a sum of
shifts of elements of $\Pcb_\Vb$ we have also
$\Yb_\Lc\in\Rb_\nu$-$\projb$ for each $\Lc\in\Qcb_\Vb$. In other
words, we have constructed an additive functor
$$\Yb_\Vb:\Qcb_\Vb\to\Rb_\nu\text{-}\projb,\quad
\Lc\mapsto\Yb_\Vb(\Lc)=\Yb_\Lc$$ which commutes with the shift of
the grading. Let $\Yb$ denote the $\Ac$-linear map
$$\Yb:K(\Qcb)\to K(\Rb),\quad
[\Lc]\mapsto[\Yb_\Lc].$$ Since $\Bc_\nu$ is a basis of the
$\Ac$-module ${}_\Ac\fb_\nu$, there is a unique $\Ac$-linear map
$${}_\Ac\fb_\nu\to K(\Rb_\nu),\quad
\bb_\Lc\mapsto[\Yb_\Lc],\quad\forall\Lc\in\Pcb_\Vb.$$ Taking the
limit over all $\Vb\in\Vcb$ we get a $\Ac$-linear map
$$\mu_\Ac:{}_\Ac\fb\to K(\Rb)\leqno(4.9)$$
such that $\mu_\Ac\circ\l_\Ac=\Yb$. We claim that we have
$$\mu_\Ac=\g_\Ac.$$
Indeed (4.4), (4.8) yield
$$\aligned
\mu_\Ac(\theta_{\ib})= \mu_\Ac\l_\Ac([{}^\delta\!\Lc_\ib])=[\Rb_\ib]
=\g_\Ac(\theta_{\ib}),\quad\forall\ib\in I^\nu.
\endaligned$$
By base change the $\Ac$-linear maps $\g_\Ac$, $\mu_\Ac$ yield
$\Kc$-linear maps
$$\g:\fb\to K(\Rb)_{\Kc},\quad
\mu:\fb\to K(\Rb)_{\Kc},\quad K(\Rb)_{\Kc}=K(\Rb)\otimes_\Ac\Kc.$$
Since the $\Ac$-modules ${}_\Ac\fb$, $K(\Rb)$ are both torsion free
it is enough to prove that the maps $\mu$, $\g$ are the same. This
is obvious because the monomials $\theta_{\ib}$ span the
$\Kc$-vector space $\fb_\nu$ as $\ib$ runs over $I^\nu$.

Next we prove that $\mu_\Ac$ takes the elements of $\Bc_\ZZ$ to the
classes in $K(\Rb)$ of the indecomposable projective objects. We
claim that it is enough to prove that $\Yb_\Lc$ is indecomposable in
$\Rb_\nu$-$\projb$ for each $\Lc\in\Pcb_\Vb$. Indeed the map
$\g_\Ac$ is invertible and $\mu_\Ac=\g_\Ac$. Thus $\mu_\Ac$ takes
$\Bc_\ZZ$ to a $\ZZ$-basis of $K(\Rb)$. Since
$\mu_\Ac\circ\l_\Ac=\Yb$, if $\Yb_\Lc$ is indecomposable in
$\Rb_\nu$-$\projb$ for each $\Lc\in\Pcb_\Vb$, then $\mu_\Ac$ takes
$\Bc_\ZZ$ to a subset of the $\ZZ$-basis of $K(\Rb)$ consisting of
the indecomposable projective objects. So we are done.

To prove the claim we'll prove that the top of  $\Yb_\Lc$ in
$\Rb_\nu$-$\modb$ is a simple object. Fix a simple quotient
$\Yb_\Lc\to L$ in $\Rb_\nu$-$\modb$. Note that $L$ is finite
dimensional, because $\Rb_\nu$ is finite dimensional over its
center. The center of $\Rb_\nu$ is equal to $\Sb_\nu=\Sb_\Vb$. See
\cite{KL1, sec.~2.4}. Let $\Sb_\nu^+\subset\Sb_\nu$ be the unique
graded maximal ideal. It acts by zero on each simple finite
dimensional left graded $\Rb_\nu$-module. So it acts by zero on $L$.
We set
$$\Rb_\nu^0=(\Sb_\nu/\Sb^+_\nu)\otimes_{\Sb_\nu}\Rb_\nu,\quad
\Yb_\Lc^0=(\Sb_\nu/\Sb^+_\nu)\otimes_{\Sb_\nu}\Yb_\Lc.$$ Then $L$ is
a simple left graded $\Rb_\nu^0$-module, $\Yb^0_\Lc$ is a projective
left graded $\Rb_\nu^0$-module and, since taking the tensor product
with $\Sb_\nu/\Sb_\nu^+$ is a right exact functor, the surjective
map $\Yb_\Lc\to L$ factors to a surjective left $\Rb_\nu^0$-module
homomorphism $\Yb_\Lc^0\to L.$

To simplify the notation we'll use the same symbol for a complex in
$\Dcb_{G_\Vb}(E_\Vb)$ and its image  by the forgetful functor
$\Dcb_{G_\Vb}(E_\Vb)\to \Dcb(E_\Vb)$. Write
$$\Rb^1_\nu=\Ext^*_{}({}^\delta\!\Lc_\Vb,
{}^\delta\!\Lc_\Vb),\quad
\Yb^1_\Lc=\Ext^*_{}({}^\delta\!\Lc_\Vb,\Lc).$$ The proof of the
following lemma is postponed in Section 4.10.

\proclaim{4.8. Lemma} The forgetful map $\Rb_\nu\to\Rb_\nu^1$
factors to a graded $\k$-algebra isomorphism
$\Rb_\nu^0\to\Rb_\nu^1$. The forgetful map $\Yb_\Lc\to\Yb^1_\Lc$
factors to a left graded $\Rb_\nu^0$-module isomorphism
$\Yb_\Lc^0\to\Yb^1_\Lc$.
\endproclaim

For each simple perverse sheaf $\Lc'\in\Pcb_\Vb$ we define a finite
dimensional $\CC$-vector space $M_{\Lc'}$ by
$M_{\Lc'}=\bigoplus_{r\in\ZZ}M_{\Lc',r}$ with
$${}^\delta\!\Lc_\Vb=\bigoplus_{\Lc'\in\Pcb_\Vb}\bigoplus_{r\in\ZZ}
M_{\Lc',r}\otimes\Lc'[r].$$ We have
$$\Rb_\nu^0=\Rb_\nu^1=
\bigoplus_{\Lc',\Lc''\in\Pcb_\Vb}\Hom(M_{\Lc''},M_{\Lc'})\otimes
\Ext^*(\Lc'',\Lc').$$ Thus $\Rb_\nu^0$ has a natural structure of
finite dimensional $\ZZ_{\geqslant 0}$-graded $\k$-algebra whose
degree zero part is a semi-simple $\k$-algebra isomorphic to
$\bigoplus_{\Lc'\in\Pcb_\Vb}\End(M_{\Lc'})$. See \cite{CG, sec.~8}
for more details. In particular each $\CC$-vector space $M_{\Lc'}$
has a natural structure of simple left $\Rb_\nu^0$-module. Now, we
have
$$\Yb_\Lc^0=
\bigoplus_{\Lc'\in\Pcb_\Vb}M_{\Lc'}\otimes \Ext^*(\Lc',\Lc),$$ a
$\ZZ_{\geqslant 0}$-graded $\Rb_\nu^0$-module whose degree zero part
is isomorphic to $M_\Lc$. Observe that $\Yb_\Lc^0$ is generated by
its degree zero subspace as a $\Rb_\nu^0$-module, because
$$\aligned
\Rb_\nu^0\star (\Yb_\Lc^0)_0 &=
\Bigl(\bigoplus_{\Lc',\Lc''\in\Pcb_\Vb}\Hom(M_{\Lc''},M_{\Lc'})\otimes
\Ext^*(\Lc'',\Lc')\Bigr)\star
\Bigl(M_{\Lc}\otimes\Ext^*(\Lc,\Lc)\Bigr)\cr
&=\bigoplus_{\Lc'\in\Pcb_\Vb}M_{\Lc'}\otimes
\bigl(\Ext^{*}_{}(\Lc',\Lc)\star\Ext^*(\Lc,\Lc)\bigr) \cr
&=\Yb_\Lc^0.
\endaligned$$
Thus the $\Rb_\nu^0$-module $\Yb_{\Lc}^0$ has a unique simple
quotient which is isomorphic to $M_\Lc$. Hence we have $L=M_{\Lc}$.
Thus $\Yb_\Lc$ has a unique simple quotient in $\Rb_\nu$-$\modb$. We
are done.

\qed

\vskip3mm

\subhead 4.9.~Remarks\endsubhead $(a)$ Fix a pair $\ibt=(\ib,\ab)$
in $Y_\nu$. Let $\ib'\in I^\nu$ be the sequence obtained by
expanding the pair $\ibt$. By (4.2), (4.7) and (4.8) we have
$$[\ab]![\Rb_{\ibt}]= [\ab]![\Yb_{{}^\delta\!\Lc_{\ibt}}].$$ Since the
$\Ac$-module $K(\Rb_\nu)$ is torsion free this implies that
$[\Rb_{\ibt}]=[\Yb_{{}^\delta\!\Lc_{\ibt}}]$. Now, recall that any
object in $\Rb_\nu$-$\projb$ is a direct sum of a finite number of
indecomposable ones, and that the indecomposable projective objects
yield a basis of the Abelian group $K(\Rb_\nu)$. Therefore we have
$\Rb_{\ibt}\simeq \Yb_{{}^\delta\!\Lc_{\ibt}}$ in
$\Rb_\nu$-$\projb$.

\vskip1mm

$(b)$ By (4.3), (4.6) and the previous remark we have
$$[\Yb_{\Lc\circledast\Lc'}]=
[\Ind_{\nu,\nu'}(\Yb_\Lc\otimes\Yb_{\Lc'})],\quad
\forall\Lc\in\Qcb_\Vb, \Lc'\in\Qcb_{\Vb'}.$$ Therefore, the same
argument as above yields a non canonical isomorphism in
$\Rb_{\nu+\nu'}$-$\projb$
$$\Yb_{\Lc\circledast\Lc'}\simeq\Ind_{\nu,\nu'}(\Yb_\Lc\otimes\Yb_{\Lc'}).$$

\vskip1mm

$(c)$ Note that we dont use the surjectivity of $\g_\Ac$; on the contrary this surjectivity
follows from our arguments. Namely, the image of $Y$ is the free abelian group
spanned by the subset of the basis of indecomposable objects given by $[Y_\Lc]$,
with $\Lc\in \Pcb_\Vb$ for some $\Vb.$ On the other hand the image of $Y$ is the 
$\Ac$-span of the $[\Rb_\ib]$'s, with $\ib\in I^\nu$ and $\nu\in \NN I.$
Since any indecomposable projective module is a direct summand of $\Rb_\ib$ for some 
$\ib$ (up to a shift of the grading), the image of $Y$ must contain all the indecomposable projective
objects.

\subhead 4.10\endsubhead {\sl Proof of Lemma 4.8 :} To prove the
first claim we must check that the forgetful map yields an
isomorphism
$$(\Sb_\Vb/\Sb^+_\Vb)\otimes_{\Sb_\Vb}H_*^{G_\Vb}(Z_\Vb,\k)\to
H_*(Z_\Vb,\k).$$ Let $\Pb_\Vb^+\subset\Pb_\Vb$ be the unique graded
maximal ideal. By (2.3) it is enough to prove that the forgetful map
yields an isomorphism
$$(\Pb_\Vb/\Pb^+_\Vb)\otimes_{\Pb_\Vb}H_*^{T_\Vb}(Z_\Vb,\k)\to
H_*(Z_\Vb,\k).$$ This is well-known. It is proved using a
decomposition of $Z_\Vb$ into affine cells. See \cite{CG, chap.~6}
for similar results for the equivariant K-theory of the Steinberg
variety, and \cite{GKM, sec.~1} for details on equivariantly formal
$T_\Vb$-varieties.

The second claim follows from the first one. Indeed, the forgetful
map
$$(\Sb_\Vb/\Sb_\Vb^+)\otimes_{\Sb_\Vb}
\Ext^*_{G_\Vb}({}^\delta\!\Lc_\Vb,{}^\delta\!\Lc_\Vb)\to
\Ext^*_{}({}^\delta\!\Lc_\Vb,{}^\delta\!\Lc_\Vb)$$ is invertible.
Further $\Ext^*_{G_\Vb}({}^\delta\!\Lc_\Vb,\Lc)$ is a direct summand
of $\Ext^*_{G_\Vb}({}^\delta\!\Lc_\Vb,{}^\delta\!\Lc_\Vb)$ and
$\Ext^*_{}({}^\delta\!\Lc_\Vb,\Lc)$ is a direct summand of
$\Ext^*_{}({}^\delta\!\Lc_\Vb,{}^\delta\!\Lc_\Vb)$. Therefore the
forgetful map yields also an isomorphism
$$(\Sb_\Vb/\Sb_\Vb^+)\otimes_{\Sb_\Vb}
\Ext^*_{G_\Vb}({}^\delta\!\Lc_\Vb,\Lc)\to
\Ext^*_{}({}^\delta\!\Lc_\Vb,\Lc).$$

\qed


\Refs \widestnumber\key{ABCD}

\ref\key{BBD}\by Beilinson, A., Bernstein, J., Deligne, P.\book
Faisceaux pervers\bookinfo Ast\'erisque\yr 1982\endref

\ref\key{BL}\by Bernstein, J., Lunts, V.\book Equivariant sheaves
and functors \yr 1994\bookinfo LNM\vol 1578\publ Springer\endref

\ref\key{CG}\by Chriss, N., Ginzburg, V.\book Representation theory
and complex geometry\publ Birkhauser\yr 1997\endref


\ref\key{D}\by Deligne, P.\paper Theorie de Hodge, III\jour
Publications Math\'ematiques de l'IHES\vol 44\yr 1974\pages
5-77\endref

\ref\key{GKM}\by Goresky, M., Kottwitz, R., MacPherson, M.\paper
Equivariant cohomology, Koszul duality, and the localization
theorem\jour Invent. Math.\vol 131\yr 1998\pages 25-83\endref

\ref\key{J}\by Joshua, R\paper modules over convolution algebras
from equivariant derived categories\jour J. of Alg.\vol 203\yr 1998
\pages 385-446
\endref

\ref\key{KL1}\by Khovanov, M., Lauda, A. D.\paper A diagrammatic
approach to categorification of quantum groups I\jour preprint
arXiv:0803.4121\endref

\ref\key{KL2}\by Khovanov, M., Lauda, A. D.\paper A diagrammatic
approach to categorification of quantum groups II\jour preprint
arXiv:0804.2080\endref


\ref\key{L1}\by Lusztig, G.\paper Quivers, perverse sheaves, and
quantized enveloping algebras\yr 1991\pages 365-421\vol 4\jour J.
Amer.Math. Soc.\endref

\ref\key{L2}\by Lusztig, G.\book Introduction to quantum groups \vol
110\publ Birkh\"auser\bookinfo Progress in Mathematics\yr
1993\endref

\ref\key{L3}\by Lusztig, G.\paper Cuspidal local systems and graded
Hecke algebras, I\jour Publications Math\'ematiques de l'IHES\yr
1988\pages 145-202\vol 67\endref

\ref\key{L4}\by Lusztig, G.\paper Cuspidal local systems and graded
Hecke algebras, II\jour Canadian Mathematical Society Conference
Proceedings \yr 1995\pages 217-275\vol 16\endref

\ref\key{R}\by Rouquier, R. \paper 2-Kac-Moody algebras \jour
preprint arXiv:0812.5023\endref





































































\endRefs
\enddocument